\documentclass[a4paper,11pt]{amsart}
\usepackage{graphicx}
\usepackage{amsmath}
\usepackage{amssymb}
\usepackage{oldgerm}
\usepackage[all]{xy}

\renewcommand{\Im}{\operatorname{Im}\nolimits}
\renewcommand{\mod}{\operatorname{mod}\nolimits}

\newcommand{\id}{\operatorname{id}\nolimits}

\newcommand{\N}{\operatorname{\mathbb{N}}\nolimits}

\newtheorem{theo}{Theorem}[section]

\newtheorem{cor}[theo]{Corollary}

\newtheorem{lemma}[theo]{Lemma}
\newtheorem{prop}[theo]{Proposition}
\newtheorem{defi}[theo]{Definition}
\newtheorem{rem}[theo]{Remark}

\theoremstyle{definition}

\theoremstyle{definition}
\newtheorem{exa}[theo]{Example}

\theoremstyle{remark}

\theoremstyle{definition}

\theoremstyle{definition}

\begin{document}
\title{Classification of pointed rank one Hopf Algebras}
\author{Sarah Scherotzke}
\address{Sarah Scherotzke \newline Mathematical Institute \\ University of Oxford \\ 24-29 St.\
Giles \\ Oxford OX1 3LB \\ United Kingdom} 

\email{scherotz@maths.ox.ac.uk}
\date{\today}

 \maketitle

\begin{abstract}
We classify pointed rank one Hopf algebras over fields of prime characteristic which are generated as algebras by the first term of the coradical filtration.
%If $H_0$ is a sub-Hopf algebra, which is the case if $H$ is a pointed Hopf algebra, then $H_1$ is a $H_0$-Hopf module and therefore free over $H_0$. 
%We say $H$ is of rank one if $H_1$ has dimension 2 as a free $H_0$-module.
We obtain three types of Hopf algebras presented by generators and relations. For Hopf algebras with semi-simple coradical only the first and second type appears. We determine the indecomposable projective
modules for certain classes of pointed rank one Hopf algebras.  
\end{abstract}

\section{Introduction}

In \cite{KR} Krop and  Radford introduce the rank as a measure of complexity
for finite-dimensional Hopf algebras whose coradical is a sub Hopf algebra. They classify the finite-dimensional pointed Hopf algebras $H$ of rank one which are generated as algebras by the first term of the coradical filtration $H_1$ over fields $k$ of characteristic zero.
These algebras are parametrized by tuples $(G, \chi, a, \alpha)$, where $G$ is a finite group, $\chi$ a linear charakter of $G$, $a\in Z(G)$ and $\alpha \in k$. As algebras they are generated by a $(a,1)$-primitive element $x$ and the
group $G$, subject to the action $gxg^{-1}=\chi(a)x$ and $x^n=\alpha(a^n-1)$. 

%A basis of $H_D$ is given by $\{x^ig|g\in G, 0\le i \le n-1\}$. For any datum $D$ and the given relations and generators there exists a uniquely determined Hopf algebra $H_D$ that is also pointed and rank one. It is easy to describe for which tuples $D$ and $D'$ the given Hopf algebras $H_D$ and $H_{D'}$ are isomorphic. 

In this paper, we classify the finite-dimensional pointed rank one Hopf
algebras $H$ over algebraically closed fields $k$ of prime characteristic $p$, which are generated as algebras by
$H_1$. The list consists of three types of Hopf algebras parametrized by tuples. The type only depends on the conjugation action of the so called skew point $a$ on the $(a,1)$-primitive elements. The first type is parametrized by the same tuple as the one in Krop and Radford's paper and is given by similar relations. For Hopf algebras with semi-simple coradical only the first and second type appears. The second and third types are new. Analyzing the third type requires new ideas and methods. In every type we analyse when two Hopf algebras parametrized by a different tuple are isomorphic as Hopf algebras. Furthermore, we are able to determine the indecomposable projective
modules for certain Hopf algebras in every type. The
second and third type provide new examples of pointed rank one Hopf algebras.

%We can view $k$ as right $H_0$ module via $1*h=\epsilon(h)1$ for all $h \in H_0$. We call $H$ a rank $n$ Hopf algebra if $\dim k\otimes_{H_0} H_1= n+1$. 
By \cite {TW} and \cite{KR} we know that a finite-dimensional pointed rank one Hopf algebra has a uniquely determined group-like element $a$ such that there exists a $(a,1)$-primitive elment which does not lie in the coradical. Let $X$ be the set of $(a,1)$-primitive elments which do not lie in the coradical. The group of group-like elements $G$ together with any $x\in X$ generate the algebra and determine its Hopf algebra structure. 

To describe the first type we assume that there exists an $x\in X$ which is an eigenvector with eigenvalue $q\not =1$ of the conjugation action of $a$. Then we can show that $x$ is a simultaneous eigenvector of the conjugation action of $G$. The eigenvalues of this action are given by a linear character $\chi$ of $G$. Similarly to \cite{KR} we compute $n$ such that $\{gx^i|g\in G, 0\le i\le n-1\}$ is a basis of the algebra. We call $n$ the degree of $x$. We show that $n$ is equal to the multiplicative order of $q$ and that $x^n=\alpha(a^n-1)$ for an $\alpha \in k$. Thus this Hopf algebra is parametrized by the tuple $R=(G,a,\chi,\alpha)$. 

The second type appears only if $a$ has an eigenvector $x\in X$ with eigenvalue $1$. Then the degree of $x$ is $p$ and $x^p=\alpha_{-1}(a^p-1)+\alpha_0x$ for $\alpha_{-1},\  \alpha_0\in k$. The conjugation action of $G$ on $x$ is described by two functions, a linear character $\chi$ of $G$ and a map $c:G\to k$ such that $gxg^{-1}=\chi(g)x+c(g)(a-1)$. 
Thus this Hopf algebra is parametrized by the tuple $F=(G,a,\chi,c, (\alpha_{-1}, \alpha_0))$. 

Generalyzing the results in \cite{KR}, we show that if the order of the skew point is not divisible by $p$ then there exists such an eigenvector $x\in X$ and the Hopf algebra is of first or second type. A corollary is that if $H$ has a semi-simple coradical then it is of first or second type. If $H$ is cocommutative it is of second type.

The third type appears if no eigenvector of the conjugation action of
$a$ in $X$ exists. Then we can show that there is an $x\in X$ such that $axa^{-1}=x+(a-1)$. We show that the degree of $x$ is $p$, that $x^p=x$ and that the conjugation action of $G$ is described by a group homomorphism $c:G\to (k,+)$ such that $gxg^{-1}=x+c(g)(a-1)$. Thus this Hopf algebra is parametrized by the tuple $E=(G,a,c)$. 

Conversely for every tuple $R$, $F$ and $E$ we show that there exists a uniquely determined Hopf algebra $H_R$, $H_F$ and $H_E$ which satisfies the given relations. Then we prove that $H_R$, $H_F$ and $H_E$ are pointed rank one Hopf algebras and determine for which given tuples they are isomorphic. 
 
Following \cite{KR} we define $H_R$ and $H_F$ to be of nilpotent type if $x^n=0$. We describe the projective
idecomposable modules and the simple modules for Hopf
algebras of nilpotent type. We show that $H_R$ and $H_F$ are uniserial if their coradical is
semi-simple generalyzing a result in \cite{KR} formulated in the case that the group of group-like elements is abelian and the charactersitic of the field is zero. Finally we prove that the blocks of Hopf algebras $H_R$ that are not of nilpotent type with semi-simple coradical are Morita equivalent to the Taft algebra $\Gamma_{N, \chi(a)}$ or the matrix ring $M_N(k)$. 

\bigskip

I would like to thank my supervisor Karin Erdmann for all her help on
earlier drafts of this article. I would also like to thank the referee for
her/his suggestions and comments which were extremely useful in
completing the article.     

\section{Preliminaries}
We first introduce some basic definitions and facts that will be used in
the other sections. Throughout this section $H$ is a finite-dimensional Hopf algebra.
\begin{defi}\cite[Definitions 5.1.5, 5.2.1]{M}
The {\rm coradical} $H_0$ of $H$ is the sum of its
simple subcoalgebras. We define inductively $H_i:=\Delta^{-1}(H\otimes H_{i-1} +H_0 \otimes H)$ and call $(H_i)_{i\in  \N}$ the {\rm coradical filtration} of $H$. 
A coradical filtration is called a {\rm Hopf algebra filtration} if $H_iH_j \subset H_{i+j}$ 
and $S(H_i) \subset H_i$ for all $i,j\in \N$.
\end{defi}
The following results are well known.
\begin{lemma}\cite[Theorem 5.2.2, Lemma 5.2.8]{M}
Let $(H_i)_{i\in  \N}$ be the {\rm coradical filtration} of $H$. 
\begin{enumerate}
 \item The $(H_i)_{i\in  \N}$ are a family of
subcoalgebras satisfying $H_i\subset H_{i+1}$for all $i\in \N$ and $\bigcup_{i\in \N}
H_i=H$. For the comutiplication of the elements of the filtration we have $\Delta(H_i)\subset \sum_{k=0}^i H_k\otimes H_{i-k}$. 
\item The coradical filtration $(H_i)_{i\in \N}$ is a Hopf algebra filtration if and only if $H_0$ is a sub Hopf algebra of $H$. 
\end{enumerate}
\end{lemma}
If $H$ is a pointed Hopf algebra, then $H_0$ is the group algebra generated by the set of group-like elements, which is a sub Hopf algebra of $H$. Thus the 
coradical filtration of $H$ is a Hopf filtration. 

We introduce variations of primitive elements of a Hopf algebra. We will use certain primitive elements in order to classify pointed rank one Hopf algebras. 
\begin{defi} \cite[Definition 1.3.4]{M}
Let $c\in H$ be an element such that $\Delta
(c)=c\otimes g+h \otimes c$ for some group-like elements $g,h \in G(H)$.
Then we say that $c$ is a {\rm $(g,h)$-primitive element}, and denote the set of
those elements by $P_{g,h}$. If $c\in P_{g,1}$,
then we say that $c$ is a {\rm skew primitive element}.
\end{defi}In the following lemma we determine the skew primitive
elements in $kG(H)$.  
\begin{lemma}\label{skew} Let $H$ be a Hopf algebra, $a\in G(H)$ a group-like element and $c\in kG(H)$ a skew primitive element in $P_{a,1}$. Then $c=\alpha(a-1)$ for some $\alpha \in k$. Thus $P_{a,1}\cap kG(H)=\{\alpha(a-1)|\alpha \in k\}.$
\end{lemma}\begin{proof}

The group-like elements are linearly independent by \cite[3.2.1]{S}. Thus we can write $c$
as a uniquely determined linear combination of the group-like elements.
We take $c=\sum_{g\in G(H)}b_gg$. Then the following equation holds
\[\sum_{g \in G(H)}b_g(g\otimes a+1\otimes g)= \Delta(c)=\sum_{g\in
G(H)}b_gg\otimes g.\]As the elements of $G(H)\otimes_k G(H)$ are a basis
of $kG(H)\otimes_k kG(H)\cong k(G(H)\otimes_k G(H))$, comparing coefficients of this basis on the left and right hand side gives us $b_g=0$ for all $g\in G(H)-\{a,1\}$ and $b_a=-b_1$.\end{proof}
We define certain kinds of group-like elements.

\begin{defi}We say that $g \in G(H)$ is a {\rm  skew point} if $P_{g,1} \setminus H_0$ is not empty.
\end{defi}

The following lemma shows that every pointed, non cosemi-simple Hopf algebra $H$ has a skew point $a\in G(H)$. This result can also be found in a similar form in \cite[Lemma 1, Prop. 1]{KR}.

\begin{lemma}\label{H_1}
Let $H$ be a pointed  Hopf algebra %coalgebra,
and $(H_i)_{i\in \N} $ its Hopf filtration. The set of primitive elements lies in $H_1$. Furthermore we have $H_1=H_0+\sum_{\{a,b\in G(H)\}} P_{a,b}.$ So if $H\not =H_0$ there exists a skew primitive element $x \in
P_{a,1} \setminus H_0$ for some $a\in G(H)$.
\end{lemma}
\begin{proof} As all group-like elements lie in
$H_0$, the primitive elements are in $H_1$ by definition. The second part is \cite[Prop. 2]{TW}.
\end{proof}
We introduce the rank of a pointed Hopf algebra and some properties of
these Hopf algebras.  

\bigskip 

Let $H$ be a pointed Hopf algebra and $(H_i)_{i\in \N}$ the coradical filtration of $H$. Then $(H_i)_{i\in \N}$ is a
Hopf filtration and we can view the $H_i$'s as $H_0$-module. The trivial module $k$ is a $H_0$-module via the counit map. If $\dim k \otimes _{H_0} H_1 = t+1$ we call $H$ a {\rm rank} $t$ pointed Hopf algebra. Using the next remark, we see that $H_i$ is a free left $H_0$-module for all $i\in \N$.
 
\begin{rem}\cite[section 1]{KR}\label{free}
Let $V_0\subset V_1\subset \ldots $ be a chain of subcomodules of $H$ such that $V_0$
is a sub Hopf algebra, $V_0V_i\subset V_i$ and $\Delta(V_i)\subset \sum_{t=0}^{i}V_t
\otimes V_{i-t}$. Then we can view $Q_i:= V_i/V_{i-1}$ as a left $V_0$-Hopf
module via the map $\rho_i:Q_i\to V_0\otimes Q_i, h+V_{i-1} \mapsto
\Delta (h)+ V_i\otimes V_{i-1}$. By the Fundamental Theorem of Hopf
modules \cite[Theorem 4.1.1]{S} either $Q_i=\{0\}$, or it is 
a free left $V_0$-module on a linear basis of ${Q_i}^{coV_0}=\{z\in Q_i|\rho_i(z)=1\otimes z\}$. Therefore $V_i$ is a free
left $V_0$-module for all $i\in \N$. 
\end{rem}

Let $a$ be a skew point of a pointed rank one Hopf algebra and $x \in
P_{a,1}\setminus H_0$, then ${1,x}$ is a basis of $H_1$ as left $H_0$-module by the
preceding remark. 

\bigskip

The next lemma shows that if $H$ is a pointed, rank one Hopf algebra
then it has exactly one skew point. This skew point is an element of the center of $G(H)$. The proof of the following is a combination of \cite[Lemma 1(c)]{KR} and \cite[Prop. 1(d)]{KR}. However the setup is slightly different and therefore we give the details. 
\begin{lemma}\label{xa}
Let $H$ be a pointed, rank one Hopf algebra and
$x \in P_{a,1} \setminus H_0$ for some $a\in G(H)$. Then the set of skew primitive elements in $H_1-H_0$ 
is given by $P_{a,1} \setminus H_0:=\{sx+t(a-1)|s\in k^*, t
\in k\}.$ Thus $a\in G(H) $ is uniquely determined.
We have $a\in Z(G)$
and for each $y \in P_{a,1} \setminus H_0$ there are maps $\chi$ and $c$ depending on $y$ such that
$gyg^{-1}=\chi(g)y+c(g)(a-1)$ for all $g\in G$. Then $\chi:G \to k^*$ is
a linear character and $c:G\to k$ satisfies the condition
$c(hg)=\chi(g)c(h)+c(g)$. 
\end{lemma}
\begin{proof} 
Suppose that $z \in H - H_0$ is a skew primitive element in $P_{g,1}$ for some $g\in G$. Then $z\in H_1 -H_0$ and 
has a unique presentation as 
 $z=a_0 +a_1x$ for $a_0,a_1\in H_0$ and $a_1\not = 0$. Thus $\Delta(a_1)(x\otimes
a)+\Delta(a_1)(1\otimes x) +\Delta(a_0)=\Delta(z)=z \otimes g +1\otimes
z=a_0\otimes g+1\otimes a_0+(a_1\otimes g)(x\otimes 1)+(1\otimes a_1)(1\otimes
x)$. Using the fact that $\{1\otimes1, x\otimes 1, 1\otimes x, x \otimes x \}$ is a basis of 
 $H_1\otimes H_1$ as a $H_0\otimes H_0$-module gives us $g=a$ and
therefore $\Delta(a_0)=a_0\otimes a+1\otimes a_0$ and
$\Delta(a_1)=1\otimes a_1$. Thus by applying Lemma \ref{skew} we have
$a_0=\alpha(a-1)$ with $\alpha \in k$. Applying the counit to the last
equation gives us $a_1=1\epsilon(a_1)$. Thus $a_1\in k^*$. 

Let $y \in P_{a,1} \setminus H_0$. The element $gyg^{-1}$ is a
$\{gag^{-1},1\}$-skew primitive element and lies in $H_1- H_0$. The
previous part gives us $gag^{-1}=a$ for all $g\in G$ and
$gyg^{-1}=\chi(g)y+c(g)(a-1)$ for $\chi(g)\in k^*$ and $c(g)\in k$. The equation $\chi(hg)y+c(hg)(a-1)=hgy g^{-1}h^{-1}=\chi(h)\chi(g)y+(\chi(g)c(h)+c(g))(a-1)$ proves the relations for $\chi$ and $c$. \end{proof}

Finally we state a combinatorial result that we will use to prove some results of the next section. 
\begin{lemma}\cite[Corollarly 2]{R1}\label{q}
Let $k$ be a field of positive characteristic $p$, $q \in k$ and $n\in \N$. We set $(n)_q=\sum_{i=0}^{n-1}q^i$ and define \[\binom{n}{m}_q:=(n)_q!/(n-m)_q!(m)_q!.\]
Suppose $n>1$. Then $\binom{n}{m}_q=0$ for all $1\le m \le n-1$ if and only if 
%the characteristik of $k$ is zero and $q$ is a primitive $n$-th root of unity or the characteristik $p$ of $k$ is positive and 
$n=Np^r$ with $1 \le N$, $0 \leq r$ where $N$ and $p$ are coprime and $q$ is a primitive $N$-th root of unity.
\end{lemma}

\section{ Pointed rank one Hopf algebras of first and second type}
In this section we introduce the pointed Hopf algebras $H_R$ and $H_F$ given by tuples $R$ and $F$ and generators and relations similar to \cite{KR}. We determine under which conditions a pointed, rank one Hopf algebra is isomorphic to $H_R$ or $H_F$. In particular we classify pointed rank one Hopf algebras with semi-simple coradical and cocommutative pointed rank  one Hopf algebras.
Throughout this section we will assume that $H$ is a pointed,
finite-dimensional Hopf algebra over a field $k$ of
positive characteristic $p$ such that $H\not = H_0$. Let $a$ a the skew
point of $H$. 

\smallskip

As in characteristic zero \cite[Prop.1]{KR} we can find for any skew point $a$ whose order is not divisible by $p$, a skew primitive element in $P_{a,1} \setminus H_0$ which is an eigenvector of the
conjugation action of $a $.

\begin{lemma}\label{a}
Let $a$ be a skew point of $H$, $k$ an algebraically closed field and suppose that $p$ does not divide the
order of $a$. Then there is an element $x \in P_{a,1} -H_0$ such that
$axa^{-1}=qx$ for a primitive  $N$-th root of unity $q\in k^*$.
 If $q\not = 1$ and $H$ is rank one, then $x$ is uniquely determined up to scalar multiplicities by these properties. 
 \end{lemma}
\begin{proof}

By definition $P_{a,1}\subset H_1$ contains an element of $H_1-H_0$ and
is therefore not trivial. Let $m$ be the order of $a$. For the linear
transformation $T:P_{a,1} \to P_{a,1}, y \mapsto aya^{-1}$ we have then
$T^m=\id$. As $p$ does not divide $m$ and as $k$ is algebraically
closed, $T$ is diagonalisable. Thus $P_{a,1}$ has a basis of
eigenvectors for $T$ whose eigenvalues are roots of unity. As $P_{a,1}\setminus H_0$ is
non-trivial, there is an eigenvector $x\in H_1 - H_0$ with eigenvalue
$q$. If $H$ has rank one and $q\not =1$ the uniqueness of $x$ follows
from Lemma \ref{xa}.
 
\end{proof}
We define the degree of an element in $H_1-H_0$ as the freeness degree of its powers over $H_0$. 
\begin{defi}We call the {\rm degree} of $x\in H_1-H_0$
the smallest integer $n$ such that $x^n\in \sum_{i=0}^{n-1} H_0x^i$.
Such an $n$ exists as $H$ is finite-dimensional. 
If $H$ is rank one, then all elements in
$H_1-H_0$ have the same degree. This holds because $H_1=H_0+H_0y$ for any element $y\in H_1-H_0$ and
therefore the dimension of the subalgebra generated by $H_1$ is the product of the dimension of $H_0$ and the degree of $y$. 
\end{defi} 
Suppose that $H$ is a rank one Hopf algebra which is generated as an
algebra by $H_1$ and has an element $x$ defined as
in Lemma \ref{a}. Then by \cite[Lemma 1]{KR} the degree of $x$ is the
smallest integer such that $1,x, \ldots , x^{n-1}$ form a basis for $H$
over $H_0$. 

We show that $x$ is a common eigenvector of the conjugation action of the group-like elements if $q \not =1$. 
The analogous result for a field of characteristic zero can be found in \cite[Prop. 2]{KR}. Part $(1)$ and $(2)$ for $q\not =1$ 
can be found in \cite[pp.696-699]{R1} and part $(3)$ has the same proof as \cite[Prop. 1(d)]{KR}.
\begin{theo}\label{x}
Let $a\in G(H)$ be a skew point of $H$. Suppose there is an $x\in
P_{a,1} \setminus H_0$, which is an eigenvector
of the conjugation action of $a$ with eigenvalue $q$, where $q$ is a
primitive $N$-th root of unity. 
\begin{enumerate}
%\item \cite[Prop. 2]{KR} If the characteristic of $k$ is $0$, then $x^N=\alpha(a^N-1)$ for
%some $\alpha \in k$ and $x$ has dimension $N$. Further more $q \not =1$, so that $a$ and $x$ are uniquely determined by the 
%properties described above and $H$ is not cocommutative.
 \item If $x$ has degree $n$,
then $n=Np^r$ for some $r\in \N$ 
 and \[x^n=\alpha_{-1}(a^n-1)+\sum_{i=0}^{r-1}\alpha_ix^{Np^i}\] with
$\alpha_i \in k$. In addition if $a^{N(p^{r}-p^i)}\not = 1$ for $0\le i$, then $\alpha_i=0$.
\item If $H$ is rank one we have ($r=0$ and $N>1$)
or ($r=1$ and $N=1$) for $r$ and $N$ defined as in $(1)$. 
 \item Suppose that $H$ is rank one with $N>1$ and $r=0$. Then $x$ is a common eigenvector for the conjugation action of the elements in $G(H)$, 
that is $gxg^{-1}=\chi (g)x$ for all $g\in G$ where $\chi:G(H) \to k^*$ is a linear character. 

\item Suppose that $H$ has rank one with $N=1$ and $r=1$. All elements
of $P_{a,1}\setminus H_0$ are eigenvectors of the conjugation action of
$a$ with eigenvalue 1. If $\alpha_0=0$ and $k$ is algebraically closed there exist an element $y \in P_{a,1}\setminus H_0$ with $y^p=0$. If in addition $a $ is not of order $p$, then $y$ is a common eigenvector of the conjugation action of $G(H)$. 
\item Let $H$ be rank one and generated by $H_1$ as an algebra and $A$ be the sub
Hopf algebra of $H$ generated by $a$ and $x\in P_{a,1} \setminus H_0$, where $x$
is an eigenvector of the conjugation action of $a$. Then $A$ is a $kG$-module via conjugation. We denote the smash
product between $A$ and $kG$ by $A*kG$. The map $A*kG \to H, (h,g) \mapsto hg$ is a surjective Hopf algebra homomorphism whose kernel is generated by $ a*1-1*a $. 
 \end{enumerate}
 \end{theo}
 \begin{proof}
Let $n$ be the degree of $x$. Then $x^n$ has a unique
presentation as $x^n=\sum_{i=0}^{n-1}b_ix^i$ with $b_i \in H_0$.
 Applying $\Delta$ to both sides of the equation and substituting $x^n$ gives us \begin{eqnarray*} \sum_{i=1}^{n-1} \binom{n}{i}_{q^{-1}} x^i\otimes a^ix^{n-i}&+& 
 \sum_{i=0}^{n-1}( b_ix^i\otimes
 a^n+ 1 \otimes b_ix^i)\\&=&(\Delta(x))^n= \Delta(x^n)=\sum_{i=0}^{n-1}\Delta(b_ix^i)\\
 &=&\sum_{i=0}^{n-1}\Delta(b_i)\sum_{k=0}^i\binom{i}{k}_{q^{-1}} x^k \otimes a^kx^{i-k}.
 \end{eqnarray*}
The set $\{ x^i\otimes x^j : 0\leq i, j\leq n-1\}$ 
is a free $H_0\otimes H_0$-basis of a  $H_0\otimes H_0$-submodule of
$H\otimes H$. So we can compare coefficients. 
\begin{enumerate} 

\item Equating the coefficents of $x^i\otimes x^{n-i}$ gives $\binom{n}{i}_{q^{-1}} =0$ for $1\le i \le n-1$, 

\item equating the coefficents of $x^i\otimes 1$ and $1\otimes x^i$ for $1\leq i\leq n-1$ gives $\Delta(b_i)=b_i\otimes a^{n-i}=1\otimes b_i$,

\item equating the coefficients of $x^k\otimes x^{i-k}$ for $1<i<n$ gives $\Delta(b_i)\binom{i}{k}_{q^{-1}}=0$ for all $1 \leq k\le i-1$,

\item comparing the coefficients of $1\otimes 1$ gives $\Delta(b_0)= b_0\otimes a^n+1\otimes b_0$.
\end{enumerate}

By Lemma \ref{skew} and (4) we have that $b_0=\alpha_{-1}(a^n-1)$ for an $\alpha_{-1} \in k$. 

%Using Lemma \ref{q} and the first equation we deduce $n=N$. By the third
%equation and the injectivity of $\Delta$ we have $b_i=0$
%for all $1\leq i$. With those results the second equation is also satisfied. 

%As $H_1\not = H_0$ we have $2\le N$ and therefore $q\not =1$. But $N$
%divides the order of $a$. Thus $a$ is not the identity which is equivalent to the fact that $H$ is not cocommutative.
 
Lemma \ref{q} and (1) gives us $n=Np^r$. By (3) and Lemma \ref{q} we have $\Delta(b_i)=0$ in the case that $1<i$ and there is no $s<r $ such that $i=Np^s$. As $\Delta $ is injective we have $b_i=0$. We deduce from (2) that $\epsilon(b_i)=\epsilon
(b_i)a^{n-i}=b_i\in k$ for all $1\le i\le n-1$. In case $b_{Np^s} \not =
0$ we thus have
$a^{N(p^{r}-p^s)}=1$. Suppose now $b_1\not =0$, then $a^{n-1}=a^{Np^r-1}=1$. As $N$ divides the order of $a$ it divides $Np^r-1$. Thus we have $N=1$ in this case. This proves the degree formula. 

From now on let $H$ be rank one. Suppose ($r>1$) or ($r=1$ and $N>1$), then we have for $z= x^{Np}
$ in the first case that $\Delta(z)= z\otimes a^{Np}+1\otimes z$ and for $z=x^N$
in the second case that $\Delta(z)=z\otimes a^N+1\otimes z$ using Lemma \ref{q}. Thus $z \in H_1-H_0$ and $H$ is not a rank one Hopf algebra as ${1,x,z}$ 
are $H_0$-independent. Therefore $x$ has degree $N$ for $N>1$ or degree $p$ if $N=1$.

To prove part $(3)$ we have by Lemma \ref{xa} that $gxg^{-1}=\chi(g) x +
c(g)(a-1)$ for some maps $c, \chi:G \to k$. Furthermore we have
$gag^{-1}=a$ as $a\in Z(G)$. Conjugating with
$a$ gives
\[q\chi(g)x+qc(g)(a-1)=qgxg^{-1}=agxg^{-1}a^{-1}=q\chi(g)x+c(g)(a-1).\]
If $q\not =1$ this gives us $c(g)(a-1)=0$ for all $g\in G$ and therefore
$c(g)=0$. 

If we have $r=1$ and $N=1$, then $q=1$. Let $\alpha_0=0$ and suppose $\alpha_{-1} \not = 0$. As $k$ is algebraically closed, we have $\alpha_{-1}^{1/p}\in k$. Then $y:= x -\alpha_{-1}^{1/p}(a-1)$ is an element in $P_{a,1} \setminus H_0$
and an eigenvector of the conjugation action of $a$. It satisfies
$y^p=0$. Therefore \[0=y^p= gy^pg^{-1}=(gyg^{-1})^p=c(g)^p(a^p-1).\] If we assume $a^p\not =1$ or $a=1$, then $c(g)(a-1)=0$ for all $g \in G$.

For the last part the surjectivity follows from the fact that $x$ and $H_0$ generate $H$. By direct compuation using the basis 
$(a^lx^i)_{\{0\le i \leq n-1, 1\le l \le |a|\}}$ of $A$ and the basis $(g)_{g\in G}$ of $kG$ we can see that the kernel is generated by $(a* 1-1* a)$ and that this map is a Hopf algebra morphism.  \end{proof}

We introduce Hopf algebras given by a tuple and generators and
relations. We use them to classify the pointed rank one Hopf algebras of Theorem \ref{x} with $N>1$ and $r=0$.  
\begin{defi}[first type] \label{H_R} Let $R=(G,\chi,a, \alpha)$ be a tuple where $G$ is a finite group, $a\in Z(G)$, $\chi:G \to k^*$ is a linear character with $\chi(a)\not = 1$ and $\alpha:=(\alpha_{-1}, \alpha_0, \cdots , \alpha_{r-1})\in k^{r+1}$. Let $\chi(a)$ be a primitive $N$-th root and $n:=p^rN$.

 The Hopf algebra $H_R$ is the algebra with basis $(gx^i)_{\{g\in G,0\leq i\leq n-1\}}$ and relations $gxg^{-1}=\chi(g)x$ for all $g\in G$ and $x^n=\alpha_{-1}(a^n-1)+\sum_{i=0}^{r-1}\alpha_ix^{p^iN}$. The Hopf algebra structure is determined by $\Delta(gx^i)=(g\otimes g)(x\otimes a+1\otimes x)^i$,
$\epsilon(gx^i)=\delta_{i,0}$ and $S(gx^i)=(-xa^{-1})^ig^{-1}$ 
for all $g\in G$ and $0 \le i \le n-1$. 

We assume that $\alpha_{-1}=0$ if $a^n-1=0$. In addition we assume that  
\begin{enumerate}
\item if $\alpha_i \not =0$ for some $i \ge 0$ then $a^{N(p^{r}-p^i)}= 1$ and $\chi(g)^{N(p^{r}-p^i)}=1$ for all $g \in G$;
\item if $\alpha_{-1} \not = 0$, then $\chi(g)^{Np^r}=1$ for all $g\in G$.
\end{enumerate}
 \end{defi} 
The Hopf algebra $H_R$ is uniquely determined by the given relations. In order to show that for any tuple $R$, satisfying the conditions of Definition \ref{H_R}, the Hopf algebra $H_R$ is well defined, we give a construction of $H_R$ in the following lemma. 
\begin{lemma}[constructing $H_R$]\label{const}
 Let $(G,\chi,a,\alpha)$ be a tuple with the properties described in Definition \ref{H_R}. We set $q:= \chi(a)$, $N:=|\chi(a)|$ and $m:= |a|$ and let $A:=A_{(q,\alpha,N,m) }$ be the Hopf algebra defined in \cite[(C.6)]{R1} which is generated by the group-like element $a$ and the skew primitive element $x$. Then $A$ is a $kG$-module algebra by conjugation with the relation $gxg^{-1}=\chi(g)x$ for all $g\in G$. Let $A*kG$ be the smash product of $A$ and $kG$ whose coalgebra structure is given by the tensor coalgebra $A \otimes kG$. Then the quotient $A*kG/(1*a-a*1)$ fulfils all conditions for $H_R$.
 \end{lemma} We also analyse for which tuples $R$ and $R'$ the corresponding Hopf algebras $H_R$ and $H_{R'}$ are isomorphic. The arguments are similar to \cite[Theorem 1(c)]{KR}.
\begin{lemma}\label{H_R iso}
Let $R=(G,\chi,a,\alpha) $ and $R'=(G',\chi',a',\alpha') $ be two tuples
as in \ref{H_R}.  Then $H_R$ and $H_{R'}$ are isomorphic if and only if there is a group isomorphism $f:G \to G'$ with $f(a)=a'$, 
$\chi=\chi' \circ f$ and an element $\beta \in k^*$ such that
$\beta^{Np^i}\alpha_i=\beta^{n}\alpha_i'$ for all $0 \leq i \leq r=r'$
and $\alpha_{-1}=\beta^n\alpha_{-1}'$. 
\end{lemma}
\begin{proof}First assume there is a map $f$ and a scalar $\beta $ as
above. Let $F:H_R \to H_{R'}$ the linear map with $F(gx^i)=f(g)(\beta
x')^i$ with respect to the basis $(gx^i)_{\{g\in G,\ 0\le i \le n-1\}}$ of $H_R$. As $r=r'$ and
$|\chi(a)|=|\chi'(a')|$ we have $n=n'$ such that $F$ is obviously a vector space isomorphism and direct computation shows that it is a Hopf algebra isomorphism. 

Now suppose $F: H_R \to H_{R'}$ is a Hopf algebra isomorphism. Then
$N=N'$ and $r=r'$ as $H_R$ and $H_{R'}$ need to have the same dimension. As
group-like elements are mapped to group like elements, $f:=F|_{G}$ is a
group isomorphism from $G$ to $G'$. The map $F$ also maps skew primitive
elements to skew primitive elements. Therefore $F(x)$ is a skew
primitive element with skew point $f(a)$. Direct computation shows that
the skew primitive elements of $H_{R'}-kG'$ are the elements $a'$ with corresponding
skew primitive elements $\{\beta x'+\gamma(a'-1)|\beta \in k^*, \gamma \in k \}$ and
$a'^{Np^s}$ with corresponding skew primitive elements
\[\{ \sum_{i=0}^{r-1}\beta_i x'^{Np^i}+\gamma(a'^{Np^s}-1) |\beta_i, \gamma \in k,  \beta_i=0 \mbox{ if } a^{Np^i}\not = a^{Np^s}, \exists \beta_i\not =0 \} \]  for
$0\le s\le r-1 $. As $f(a)F(x)f(a)^{-1}=\chi(a)f(x)$ where $\chi(a)$ is a
primitive $N$-th root of unity unequal $1$, we have $f(a)=a'$. Then
$F(x)=\beta x'$ for some $\beta \in k^*$. As
$F(gxg^{-1})=f(g)\beta x'f^{-1}(g)$, we have $\chi=\chi' \circ f$. Finally
$F(x)^p=F(x^p)$ gives $\beta^{Np^i}\alpha_i=\beta^{n}\alpha_i'$ for all $0 \leq i \leq r=r'$ and $\alpha_{-1}=\beta^n\alpha_{-1}'$. 
\end{proof}
We need to introduce a second type of Hopf algebras in order to classify the case $N=1$, $r=1$ from Theorem \ref{x}. First we give a typical example of a Hopf algebra of second type which will be defined in \ref{second}. 

\begin{exa}There is a Hopf algebra $L$ over a field of characteristic $2$ generated by an elementary abelian group $G$ of order $4$ and an element $x$. Let $a,b \in G$ be the generators of $G$ such that $1,a,b, ab, ax,bx,x,abx$ is a basis of $H$. The multiplication is defined by the relations $axa=x,\ bxb=x+a-1$ and $x^2=0$. The coalgebra structure is determined by $\Delta(x)=x\otimes a+1\otimes x,\ \Delta(a)=a\otimes a,\ \Delta(b)=b\otimes b$. This gives us an example for a Hopf algebra with $\chi(a)=\chi(b)=1$ and $c(b)=1, c(a)=0$ where $c$ and $\chi$ are defined as in \ref{xa}. Thus all skew primitive elements which do not lie in $L_0$ are eigenvectors of the conjugation action of $a$ but not of $b$. 
\end{exa}
In the following definition we describe those pointed rank one Hopf algebras in a formal way and characterize them. 
\begin{defi}[second type] \label{second} Let $F=(G,\chi,c,a, (\alpha_{-1},\alpha_0))$ be a tuple, where $G$ is
a finite group, $a\in Z(G)$, $\chi:G \to k^*$ is a
linear character with $\chi(a)=1$ and $c:G\to k$ is a map with $c(hg)=\chi(g)c(h)+c(g)$ for all $g,h \in G$ and $c(a)=0$.
The Hopf algebra $H_F$ is the algebra with basis $\{gx^i|g\in G, 0\le i \le p-1\}$ and relations $x^p=\alpha_{-1}(a^p-1)+\alpha_0 x$ where $\alpha_{-1}, \alpha_0 \in k$ and $gxg^{-1}=\chi(g)x+c(g)(a-1)$ for all $g \in G$.

The Hopf algebra structure is determined by $\Delta(gx^i)=(g\otimes g)(x\otimes a+1\otimes x)^i$,
$\epsilon(gx^i)=\delta_{i,0}$ and $S(gx^i)=(-xa^{-1})^ig^{-1}$ 
for all $g\in G$ and $0 \le i \le p-1$. 

We require that $\alpha_{0}\in \{0,1\}$ and that if $\alpha_0= 0$ then $\alpha_{-1}=0$. In addition we assume that 
\begin{enumerate}
\item if $\alpha_0= 0$, then $c(g)(a^p-1)=0$ for all $g\in G$;
\item if $\alpha_{0}=1$, then $a^p=a$, $\chi(g)\in F_p^*$ and
$(\alpha_{-1}(\chi(g)-1)+c(g)^p-c(g))(a-1)=0$ for all $g\in G$. 
\end{enumerate}
\end{defi} 

The following explains the two conditions $\alpha_{0}\in \{0,1\}$ and if $\alpha_0= 0$ then $\alpha_{-1}=0$. 

\medskip

Assume $H$ is rank one and $a$ has an eigenvector $x$ in
$P_{a,1}\setminus H_0$ with eigenvalue 1. Then any element of
$P_{a,1}\setminus H_0$ is fixed by conjugation with $a$ and has degree $p$. If $\alpha_0=0$ and $k$ algebraically closed, then there exists an element $y\in P_{a,1}\setminus H_0$ with $y^p=0$ by part $(4)$ of Theorem \ref{x}.

If $\alpha_0 \not = 0$ then $wx$ with $w=(\alpha_0)^{-1/p-1}\in k^*$ is an element in $P_{a,1}\setminus H_0$ such that $(wx)^p=wx+(w\alpha_0^{-1}\alpha_{-1})(a^p-1)$. Therefore we can choose an eigenvector in $P_{a,1}\setminus H_0$ of $a$ with $\alpha_0=\alpha_{-1}=0$ or with $\alpha_0=1$.

\medskip

The proof of the next lemma shows that $H_F$ is well-defined if and only if the conditions $(1)$ and $(2)$ are satisfied.
\begin{lemma}\label{well def}The Hopf algebra $H_F$ is well-defined.
\end{lemma}
\begin{proof} In order to check that the algebra structure is
well-defined, we need to consider the equation
$\chi(g)^p\alpha_0x+(\chi(g)^p\alpha_{-1}+c(g)^p)(a^p-1)=(\chi(g)x+c(g)(a-1))^p=(gxg^{-1})^p=gx^pg^{-1}= g(\alpha_{-1}(a^p-1)+\alpha_0x)g^{-1}=\alpha_{-1}(a^p-1)+\alpha_0\chi(g)x+\alpha_0c(g)(a-1)$ for all $g\in G$. By the first part of Theorem \ref{x} we have $a^p=a$ in case that $\alpha_0 = 1$. The relations for $c$ and $\chi$ follow immediately by comparing coefficients. 

The maps $\Delta$ and $\epsilon$ are algebra homomorphism and provide $H_F$ with a coalgebra structure. An easy computation using the basis shows that $S$ is the antipode. 
\end{proof}
The Hopf algebras $H_F$ can be constructed similarly to Lemma \ref{const} as the smash product of $A$ and $kG$ with the $kG$ action $gxg^{-1}=\chi(g)x+c(g)(a-1)$ on $A$ modulo the ideal generated by $(a*1-1* a)$. 
 
 We show that $H_R$ and $H_F$ are pointed and determine under which conditions
$H_R$ and $H_F$ are rank one Hopf algebras. This depends only on the degree of
$x$. The following result can also be deduced from \cite[Prop.2]{R1}. 
\begin{prop}\label{rankone}
 The Hopf algebras $H_R$ and $H_F$ are pointed. The Hopf algebra $H_R$ is rank one if and only if $r=0$. The Hopf algebra $H_F$ is rank one. 
 \end{prop}
 \begin{proof}
 We use that by \cite[Corollary 1.74]{LR} every bialgebra $B$ over a
field generated by skew primitive elements is pointed. Therefore $H_R$ and $H_F$ are pointed. 
% and $B_0=kG(B)$. Therefore the Hopf algebra $A$ from \ref{const} is pointed. As the tensor coproduct of two pointed coalgebras is pointed, $H_R$ is pointed for every $R$.

By Theorem \ref{x} we have that if $H_R$ is rank one, then $r=0$. Suppose now that $r=0$ for $H=H_R$. We have $r=1$ for $H=H_F$. Then any element $z\in H$ has a unique presentation of the form
$\sum_{i=0}^{n-1}b_ix^i$ with $n=N>1$ in the first case and $n=p$ in the
second case and $b_i\in kG(H)$ for all $0\le i\le n-1$. The conditions
$\Delta(z)=z\otimes g +1\otimes z$ and $z\not \in H_0$ give the following equation 
\begin{eqnarray*} \sum_{i=1}^{n-1}\Delta(b_i)\sum_{k=0}^i\binom{i}{k}_{q^{-1}} x^k\otimes a^kx^{i-k}&=&\sum_{i=0}^{n-1}\Delta(b_i)\Delta(x)^i\\&=& \Delta(z)=z\otimes g +1\otimes z\\&=&\sum_{i=0}^{n-1} (b_i\otimes g)(x^i\otimes 1)+(1\otimes b_i)(1\otimes x^i).
\end{eqnarray*} This equation is equivalent to \begin{enumerate}
\item for all $1\le k \le i-1$ we have $\Delta(b_i)\binom{i}{k}_{q^{-1}}=0 $;
\item furthermore $\Delta(b_0)=1\otimes b_0+b_0\otimes g;$
\item and $\Delta(b_i)=(b_i\otimes g)(1\otimes a^{-i})=1\otimes b_i$ for all $1\le i\le n-1$.
\end{enumerate}
By Lemma \ref{q} for all $2\le i \le n-1$ there exists a natural number $k$ in $\{1,\ldots, i-1\}$ such that $\binom{i}{k}_{q^{-1}}$ is not zero. This forces $b_i =0 $ for all $2\le i \le n-1$ using the first equation.   
The second equation gives us $ b_0 \in H_0\cap P_{g,1}=\{\gamma(g-1)|\gamma\in k\}$ by Lemma \ref {skew} and the last equation shows $b_1\in k^*$ and $a^{-1}g=1$. 

Thus $z=\beta x+\gamma(a-1)$ for $\beta \in k^*$, $\gamma \in k$ and $g=a$. For all $h\in P_{b,d}\setminus H_0$
with $b,d \in G(H)$ we have $h=0$ for $d^{-1}b \not= a$ and else $h=d^{-1}(\beta x+\gamma(a-1))$. As $H_1= H_0+\sum_{b,d \in G(H)}P_{b,d}\subset
H_0+H_0x$ and $x\not \in H_0$, the dimension of $k\otimes_{H_0}H_1$ is two. Thus $H$ is a rank one Hopf algebra.
 \end{proof}

We classify certain pointed rank one Hopf algebras using $H_R$ and $H_F$. 
 \begin{theo}\label{first second type}
 Let $H$ be a rank one Hopf algebra that is generated as an algebra by $H_1$.
Assume there exists
an eigenvector in $P_{a,1} \setminus H_0$ of the conjugation action of $a$ with eigenvalue $q$. If $q\not =1$ then $H$ is isomorphic to a Hopf algebra $H_R$ given by a tuple $R$ with $r=0$. If $q=1$ and $k$ algebraically closed, then $H$ is isomorphic to a Hopf algebra $H_F$ given by a tuple $F$.
\end{theo}
\begin{proof} This follows directly by the Theorem \ref{x}, Lemma \ref{const}, Lemma \ref{well def} and Proposition \ref{rankone}.
\end{proof} 
We can now classify pointed rank one Hopf algebras with semi-simple
coradical. 
 \begin{cor}[semi-simple coradical]\label{pointed}
 Let $H$ be a rank one Hopf algebra over an algebraically closed field that is generated as an algebra by $H_1$ and has a semi-simple coradical.
Then $H$ is isomorphic to a Hopf algebra $H_R$ with $r=0$ or to a Hopf algebra $H_F$.
\end{cor}
\begin{proof} The coradical of a pointed Hopf algebra is semi-simple if
and only if the characteristic $p$ of $k$ does not divide the order of the
group of group-like elements. Thus $p$ does not divide the order of the skew point and \ref{a} holds. The rest
follows then directly by the Theorem \ref{first second type}.
\end{proof} 
Note that this proof only requires that the order of the skew point is
not divisible by $p$.

\begin{cor}Let $H$ be a cocommutative rank one Hopf algebra over an
algebraically closed field. Then $H$ is isomorphic to $H_F$. Furthermore $x$ is a common eigenvector of the conjugation action of $G$. 
\end{cor}
\begin{proof}If $H$ is cocommutative then $a=1$. Therefore every element of $P_{1,1}\setminus H_0$ is an eigenvector of the conjugation action of $a$ with eigenvalue $1$. By Theorem \ref{first second type} we have $H\cong H_F$. Furthermore $c(g)(a-1)=0$ for all $g\in G$ which proves the last statement. 
\end{proof}
We introduce nilpotent type Hopf algebras which have also been defined similarly by Krop and Radford \cite{KR}.
\begin{defi}[nilpotent type]\label{nil}
Let $H$ be a rank one Hopf algebra as in \ref{x}. We call $H$ {\rm of nilpotent type} if there exists an eigenvector $x$ in $P_{a,1}\setminus H_0$ of the conjugation action of $a$ such that $x^n=0$, where $n$ denotes the degree of $x$.
\end{defi}
\begin{cor}
The rank one Hopf algebra $H_R$ is of nilpotent type if and only if $\alpha_{-1}=0$.
 The Hopf algebra $H_F$ is of nilpotent type if and only if $\alpha_{0}=0$. 
\end{cor}
\begin{proof} As $H_R$ is rank one, we have $r=0$ by \ref{rankone}. Since $\chi(a) \not =1$ the eigenvector $x\in P_{a,1} \setminus H_0$ of the conjugation action of $a$ is uniquely determined up to scalar multiplicities. By Theorem \ref{x} we have $x^N=\alpha_{-1}(a^N-1)$ and by the definition of $\alpha_{-1}$ we have that $x^N=0$ if and only if $\alpha_{-1}=0$. 

Now consider $H_F$. All elements of $ P_{a,1}\setminus H_0$ are eigenvectors of the conjugation action of $a$. An easy computation shows that if $\alpha_0=1$, then there exist no element $y$ in $ P_{a,1}\setminus H_0$ such that $y^p=0$. If $\alpha_o=0$ then by definition $\alpha_{-1}=0$ and $H_F$ is of nilpotent type. 
\end{proof}

For certain orders of the skew point we have the following structure.
\begin{cor} Let $H$ be as in \ref{first second type} and let $a$ be an element of $p$-power order. Then $H$ is isomorphic to a nilpotent type Hopf algebra $H_F$ for a suitable tuple $F=(G,a, \chi,c,(0,0))$. If $a$ does not have order $p$, then $c$ is identical zero and $x\in H_F$ is therefore a common eigenvector of the conjugation action of $G$. 
\end{cor}
\begin{proof}As the order of $a$ is a $p$-power, we have $\chi(a)=1$. Then $H \cong H_F$ by Theorem \ref{first second type}. As $a^p\not =a$ we have by Lemma \ref{well def} that $\alpha_0=0$ and $H_F$ is therefore of nilpotent type. If $a^p \not= 1$ then it follows from Theorem \ref{x}(4) that $c(g)=0$ for all $g\in G$. 
\end{proof}
We have already determined precisely in \ref{H_R iso} under which conditions rank one Hopf algebras $H_R$ and $H_{R'}$ are isomorphic and we can now determine when nilpotent type Hopf algebras $H_F$ and $H_F'$ are isomorphic.
\begin{lemma}
%Suppose that $H_{R'}$ and $H_R$ are isomorphic rank one Hopf algebras, then there exist a group isomorphism $f:G \to G'$ with $f(a)=a'$, $\chi=\chi' \circ f$ and an element $\beta \in k^*$ such that $\alpha_{-1}=\beta^{N}\alpha_{-1}'$. 
Let $H_{F'}$ and $H_F$ be rank one Hopf algebras. If there exists a group isomorphism $f:G \to G'$ with $f(a)=a'$, $\chi=\chi' \circ f$, $c=c'\circ f$ and $\alpha=\alpha'$, then $H_F$ and $H_{F'}$ are isomorphic. 
Conversely suppose that $H_{F'}$ and $H_F$ are isomorphic rank one Hopf algebras of nilpotent type and $a$ does not have order $p$. Then there exists a group isomorphism $f:G \to G'$ with $f(a)=a'$ and $\chi=\chi' \circ f$. 
\end{lemma}
\begin{proof}  %The $H_R$-part is \cite[Theorem 1(c)]{KR}. Suppose $H_R$ and $H_{R'}$ are isomorphic via the Hopf morphism $\rho:H_R \to H%_{R'}$. Then $\rho(G)=G'$ as group-like elements are mapped to group-like 
%elements and $f:=\rho|_G:G \to G'$ is a group isomomorphism. As
%$\Delta(\rho(x))=\rho(x)\otimes \rho(a)+1\otimes \rho(x)$ and $\rho(x) \in {H_{%R'}}_1-{H_{R'}}_0$ we have $\rho(a)=a'$ and $\rho(x)=\beta x'+\gamma (a'-1)$ by% Lemma \ref{xa} for $\beta \in k^*$ and $\gamma \in k$. In addition $\rho(x)$ i%s an eigenvector of the conjugation action of $\rho(g)$ with 
%eigenvalue $\chi(g)$ for all $g\in G$. As $\chi(a) \not= 1$ it is therefore clear that $\gamma(a'-1)=0$ and $\chi =\ch%i'\circ f$. This proves the first part of the lemma. Comparing $x^N$ and $x'^N$% give the relations on $\beta $, $\alpha_{-1}$ and $\alpha'_{-1}$. 
We set $\rho(gx^i):=f(g)x'^i$. Then $\rho $ is an Hopf algebra isomorphism. 

Conversely suppose now that $H_F$ and $H_{F'}$ are isomorphic via the Hopf morphism $\rho:H_F \to H_{F'}$. Analogously to \ref{H_R iso} we have that $f:=\rho|_G:G \to G'$ is a group isomomorphism with $f(a)=a'$ and $\rho(x)=\beta x'+\gamma (a'-1)$ for $\beta \in k^*$ and $\gamma \in k$. We have $0=\rho(x^p)=\gamma^p (a'^p-1)$. Thus $\gamma(a'-1)=0$ and $\chi =\chi'\circ f$.  
\end{proof}
In the following we analyse the representations of the pointed rank
one Hopf algebras $H_R$ and $H_F$. We first consider the non-nilpotent type Hopf algebra $H_R$. 

\bigskip
 
Suppose $H_R$ is not of nilpotent type and $k$ is algebraically closed. We can assume without loss of generality that $\alpha=1$. Let $Z\subset G$ be the kernel of $\chi$ and let $\{e_i|1\le i \le s\}$ be the set of central primitive idempotents in $kZ$ corresponding to the
blocks of $kZ$. As $H_R$ is not of nilpotent type
$\chi^N(g)=1$ for all $g\in G$ by \ref{H_R} (2). Thus the cyclic group $\{\chi(g)|g\in G\}$ is
generated by $\chi(a)$ and therefore $G=\cup_{i=1}^N a^i Z$. We set
$B_i:=kZe_i$ for the blocks of $kZ$. By the definition of $Z$ and as
$a\in Z(G)$, the elements $e_i$ commute
with $G $ and $x$. We have $H_R=\bigoplus_{1\le i \le s } H_i$ where
$H_i:=H_Re_i$. Let $A_i$ be the subalgebra of $H_i$ generated by $\bar x:=xe_i$ and $\bar a:=ae_i$. Then $A_i$ is isomorphic to $A_{\chi(a),1, N,|ae_i|}$ as in \cite[(C.6)]{R1}. % with the relations $\bar a \bar x=\chi(a)\bar x \bar a$ and $\bar x^N=w_y-e_y$. 
 Then $ H_i \cong A_i\otimes_kB_i/\langle e_i\otimes \bar a^N-\bar a^N \otimes e_i \rangle$.  We have $Z(A_i)=k\langle\bar a^N \rangle$ and therefore $Z(A_i)\subset Z(B_i)$.
As $e_i$ is the unique primitive central idempotent in $B_i$ and $e_i\in Z(A_i)$, $e_i$ is primitive in $Z(A_i)$ as well. Therefore $Z(A_i)$ is indecomposable. So $\bar a^N$ is
either a scalar $we_i$ for some $w\in k$ if $\dim Z(A_i)=1$ or
$Z(A_i)\cong kC_{p^t}$ for some $t\in \N$. As an algebra with indecomposable center is indecomposable $A_i$ is indecomposable. 
In the case of $\bar a^N =e_i$ we have $(\bar x)^N =0$ and $\langle \bar a \rangle$ has
order $N$. Thus $A_i$ is isomophic to the Taft algebra
$\Gamma_{N,\chi(a)}$ and $H_i\cong \Gamma_{N,\chi(a)}\otimes_k B_i$. If
$\bar a^N=we_i$ with $w\not =1$, then $A_i\cong M_n(k)$ by
\cite[1.3]{R2} as $k$ is algebraically closed and $H_i \cong M_n(B_i)$.
As the tensor product of two indecomposable algebras is indecomposable,
$H_i$ is indecomposable if $\bar a^N$ is a scalar. In the case that
$Z(A_i)\cong kC_{p^t}$, we have that $\bar a^{Np^t}$ is a scalar and $\langle 1\otimes \bar a^N-\bar
a^N\otimes 1\rangle^{p^t}=0$. Therefore idempotents can be lifted from $ A_i\otimes_kB_i/\langle e_i\otimes \bar a^N-\bar a^N \otimes e_i \rangle$ to  $A_i\otimes_kB_i$. As $A_i\otimes_kB_i$ is an indecomposable algebra, $H_i$ is indecomposable as well.

%In the case that $\bar a^N$ has order $p^t$, the ideal $(\bar x) \subset A_i$ is nilpotent and $A_i/(\bar x) \cong k\langle \bar a\rangle /k\langle \bar a^N \rangle$ which is semi-simple. Therefore $(\bar x)$ is the Jacobson radical of $A_i$ and all simple modules are the simple modules of $k C_N$ on which $\bar x$ acts as zero. 

\bigskip

The proof of \cite[Prop. 4]{KR} works in the following case.  
\begin{theo} \label{nil rep} Let $H_R$ be a pointed rank one Hopf algebra over an
algebraically closed field which is not of nilpotent type and which has a semi-simple coradical. Then its blocks are Morita equivalent to the Taft algebra $\Gamma_{N,\chi(a)}$ or to the matrix ring $M_N(k)$.
\end{theo}
\begin{proof}We fix $i$ and set $B=B_i$, $e=e_i$, $\bar x=xe$, $\bar
a=ae$ and $A=A_i$. As $kG$ is semi-simple and $k$ algebraically closed
we have $B \cong M_n(k)$ for some $n \in \N$ by the Wedderburn-Artin Theorem. The element $\bar a^N$ lies in the center of $B$,
hence $\bar a^N=we$ for some $w\in k$. Therefore $He\cong A\otimes_kB$ and
$\bar x^N=(w-1)e$. If $w=1$, then $A$ is isomorphic to the Taft algebra
$\Gamma_{N,\chi(a)}$. Therefore $He\cong \Gamma_{N,\chi(a)}\otimes
M_n(k)$ is Morita equivalent to $\Gamma_{N,\chi(a)}$ by \cite[6.11]{ASS}. If $w\not =1$, then
$A\cong M_N(k)$ and $He$ is Morita equivalent to the full matrix ring. 
\end{proof}
The representations of Hopf algebras $H$ of nilpotent type can be described entirely in terms of those of the underlying group algebra $kG(H)$. We generalize a result in 
\cite[Prop. 3]{KR} which has only been proven for the nilpotent
type in characteristic zero and $G$ an abelian group.  
\begin{theo} \label{blocks}Suppose $H$ is a pointed rank one
Hopf algebra of nilpotent type so that $H\cong H_R$ for some tuple $R$ or $H \cong H_F$
for some tuple $F$. 

(a) The simple $H$-modules are the simple $kG$-modules on which $x$ acts as zero.

%for $H\cong H_R$ or $H\cong H_F$ with $c$ identically zero. If $H\cong H_F$ and $c$ is not identically zero, then the simple $H$-modules are the simple $kG/(a-1)kG$-modules on which $x$ acts as zero. 

(b) Assume that the coradical of $H$ is semi-simple. Let $L_{\lambda}$ be the simple module where the restriction to $kG$ has character $\lambda$, and let $e_{\lambda}$ be the primitive idempotent of $kG$ associated to $\lambda$. Then $He_{\lambda}$ is an indecomposable projective $H$-module, and every indecomposable projective $H$-module is isomorphic to some $He_{\lambda}$. 

(c) The module $He_{\lambda}$ is uniserial, with composition factors
\[L_{\lambda}, L_{\lambda \chi}, \  L_{\lambda \chi^2}, \ldots, L_{\lambda\chi^{n-1}} \]
Hence two projectives $He_{\lambda}$, $He_{\mu}$ are in the same block if and only if $\mu = \lambda \chi^i$ for some $i$. 
\end{theo}

\begin{proof} 
We will prove this theorem seperately for the case $H\cong H_R$ and $H
\cong H_F$. We suppose first that $H \cong H_R$.
(a) As $H_R$ is of nilpotent type, $(x)$ is a nilpotent ideal. Therefore $(x) \subset J(H_R)$. Note that the quotient $H_R/ (x)$ is 
isomorphic to $kG$. Thus all simple modules are isomorphic to simple
$kG$-modules on which $x$ acts as zero and $a$ as the identity. 

(b) and (c) We can take $L_{\lambda} = kGe_{\lambda}$. Consider the module $He_{\lambda}$, 
define subspaces \[M^i:= \oplus_{j\geq i} x^j(kGe_{\lambda})\]
for $j=0, 1, \ldots,  n-1$. 
These are $H$-submodules of $He_{\lambda}$ and $M^{i+1}\subseteq M^i$ for each $i$. 
The action of $H$ can be calculated explicitly on the quotient, and 
one sees that $M^i/M^{i+1}$  is isomorphic to the simple module with character
$\lambda \chi^i$. 
We claim that $He_{\lambda}$ is uniserial. Let $M$ be any submodule, 
then there is $i$ such that $M$ is contained in $M^i$ but not in $M^{i+1}$.
Then $M$ contains some element of the form \[\omega = \sum_{j=i}^{n-1} x^jm_j \]
with $m_j\in kGe_{\lambda}$ where $m_i\neq 0$ (note that the $m_j$ are unique). 
A standard argument using the nilpotency of $x$ shows that then $M=M^i$.

Since $He_{\lambda}$ is uniserial, it is indecomposable. The group algebra 
$kG$ is a direct sum of simple modules $kG=\oplus kG e_j$ where
the $e_j$ are orthogonal primitive idempotents and each 
$kGe_j$ is isomorphic to $kGe_{\lambda}$ for some $\lambda$. It follows that 
$H= \oplus He_j$, and hence each indecomposable projective $H$-module 
is isomorphic to $He_{\lambda}$ for some $\lambda$. 
All other statements of (b) and (c) follow for $H\cong H_R$.

\bigskip

We suppose now that $H \cong H_F$. If $c$ is identically zero, $(x)$ is
nilpotent of index $p$ and $(x) \subset J(H).$ In this case $H/(x) \cong kG$.
If $a$ does not have order $p$, then by the proof of Theorem \ref{x} part $(4)$ the map $c$ is identically zero. 

Suppose that $a$ has order $p$ and $c$ is not identically zero. Since
$a$ is a central element of order $p$, the ideal $(a-1)H$ is nilpotent
of index $p$. As $c$ is not identically zero we have $(a-1)H \subset
(x)$. Modulo $(a-1)H$ the elements $g\in G$ act on $x$ by multiplication
through $\chi(g)$. Therefore $(x)/(a-1)H$ is nilpotent of index $p$.
Thus $(x)^{p^2}=0$. This proves $(x)\subset J(H)$.  In this case the quotient $H_F/(x)$ is ismorphic to $kG/(a-1)kG$ as $(x)\cap kG=(a-1)kG$. Thus all simple modules are 
isomorphic to simple $kG$-modules on which $x$ acts as zero and $a$ as identity. 

(b) and (c): If $kG$ is semi-simple,
$a$ does not have order $p$ and $x$ is therefore a common eigenvector of the conjugation action
of $G$. The rest of the proof is the same then the proof of (b) and (c)
in the case that $H\cong H_R$.  
\end{proof}

\section{Pointed rank one Hopf algebras of third type}
Throughout this section $H$ is a finite-dimesional pointed rank one Hopf algebra over a field $k$ of positive characteristik $p$ which is generated as an algebra by $H_1$. We assume that $H\not = H_0$. Let $a$ be the skew point of $H$.

\smallskip

In this section we classify the pointed rank one Hopf algebras which do not have an eigenvector of the conjugacy action of $a$ in $P_{a,1}\setminus H_0$. 

We first analyse how the order of $a$ and its conjugation action determine the type of Hopf algebra and introduce the third type of pointed rank one Hopf algebras. 

\begin{theo}\label{types}
Let $a$ be the skew point of $H$ and $k$ algebraically closed. Then precisely one of the following
holds:

(i) There is some $x\in P_{a,1}\setminus H_0$ such that $axa^{-1} = qx$
where $q$ is a primitive $N$-root of unity. 

(ii) There is some $x\in P_{a,1}\setminus H_0$ with 
$axa^{-1} = x+ (a-1)$. 

Suppose (i) holds. If $q \not =1$ then $H$ is isomorphic to $H_R$,
with $r=0$.
If (i) holds and $q=1$, then $H$ is isomorphic to an algebra
$H_F$ as in \ref{second}.

Suppose (ii) holds, then $p$ divides the order of $a$. 
\end{theo}
\begin{proof} 

If $x$ is an eigenvector of $a$ then Theorem \ref{x} and Proposition \ref{rankone} give the result. 

Suppose now that $a$ does not have an eigenvector in $P_{a,1} \setminus H_0$. By \ref{a} the order of $a$ must be divisible by $p$. By Lemma \ref{xa} and \ref{skew} we know that the vector space of skew primitive elements $P_{a,1}$ is $\{cx+d(a-1)|c,d\in k \}$. Thus it has dimension two. As $k$ is algebraically closed, the conjugation action of $a$ on $P_{a,1}$ is represented by a two dimensional Jordan block. As no element of $P_{a,1} \setminus H_0$ is an eigenvector, any eigenvector must be of the form $d(a-1)\in H_0$ for some $d\in k^*$. As those elements are eigenvectors with 
eigenvalue one, the entries on the diagonal of the Jordan block must be one. Thus there is an element $ x \in P_{a,1}\setminus H_0$ so that $axa^{-1}=x+(a-1)$.
\end{proof}
Considering the previous theorem we can introduce the following notations to classify $H$.
\begin{defi}
Suppose there is some $x\in P_{a,1} \setminus H_0$ with $axa^{-1}=qx$.
If $q \not =1$, we call $H$ of {\rm first type}; and if $q=1$ of {\rm
second type}. If $P_{a,1} \setminus H_0$ contains an element $x$ with
$axa^{-1}=x+(a-1)$ then we call $H$ of {\rm third type}. 
\end{defi}Note that all three cases are different. If $k$ is
agebraically closed then we have by the previous theorem that $H$ is
either of first, second or third type. 
We have already completely described the first and second type by generators and relations. 

The conjugation action of the third type can be described as follows.
\begin{cor}\label{c}
Let $H$ be of third type. Then the action of
$G(H)$ on $H$ is given by $gxg^{-1}=x+c(g)(a-1)$ where $c:G\to (k,+)$ is a group homomorphism.
\end{cor}
\begin{proof}
By Lemma \ref{xa} the conjugation action of $G$ can be described via the maps $c$ and $\chi $ with $c(hg)=\chi(g)c(h)+c(g)$ for all $g,h\in G(H).$ We have $c(h)+1=c(ha)=c(ah)=\chi(h)+c(h)$ for all $h\in G$. Thus $\chi=1$ and $c(hg)=c(h)+c(g)$ for all $g,h \in G(H)$.
\end{proof}
In order to show that third type Hopf algebras do exist we give an example in characteristic 2 and 3. 
\begin{exa}[third type]\label{exa}
Consider the four dimensional Hopf algebra $A_2$ over a field of characteristik $2$ generated by $a,x$ with basis $\{1,a,ax,x\}$ given by the following relations:
$a^2=1$, $x^2=x$ and $axa=x+(a+1)$.
The following linear maps define the Hopf algebra structure on $A$:
\begin{itemize}
\item $\Delta(a)=a\otimes a,\ \Delta(x)=x\otimes a+1\otimes x,\ \Delta(ax)=ax\otimes 1+a\otimes ax$,
\item $\epsilon(a)=\epsilon(1)=1,\epsilon(x)=\epsilon(ax)=0$,
\item $S(a)=a,\ S(1)=1$, $S(x)=xa,\ S(ax)=x$.
\end{itemize}
Via direct computation we can check that $\Delta$ and $\epsilon$ fulfil the relations and are therefore algebra homomorphisms. The other axioms of a Hopf algebra can be checked on the basis elements.

This Hopf algebra is pointed rank one such that by Lemma \ref{xa} no skew primitive element is an eigenvector of the conjugation action of $a$. Thus it is not isomorphic to a Hopf algebra $H_R$ or $H_F$.  

As a left module $A_2$ is the direct sum of the indecomposable projective modules $P_1=A_2x=\langle x, ax \rangle $ and $A_2(x+1)=\langle x+1, a+ax \rangle $. The algebra $A_2$ has two simple one dimensional modules $S_1$ and $S_2$ where $a$ acts on both as the identity and $x$ acts on $S_1$ as the identity and on $S_2$ as zero. Then soc$ P_1=\langle  x+ax \rangle  \cong S_2$, $P_1/$rad$P_1 \cong S_1$ and soc$P_2 =\langle 1+ax+a+x \rangle  \cong S_1$, $P_2/$rad$P_2 \cong S_2$.
\end{exa} 
\begin{exa}
Let $A_p$ be the algebra over a field of characteristic $p$ generated by
$a$ and $x$ with basis $(x^ia^j)_{0\le i,j < p}$ and the relations $axa^{-1}=x+(a-1)$, $x^p=x$ and $a^p=1$. Then we define linear maps by:
\begin{itemize}
\item $\Delta(a)=a\otimes a, \Delta(x)=x\otimes a+1\otimes x, \Delta(x^ia^j)=\Delta(x)^i\Delta(a)^j$ for all $0\le j,i\le p-1$,
\item $\epsilon(x^ia^j)=0$ for $i>0$ and equal to $1$ else,
\item $S(a)=a, S(x)=-xa^{-1}$ and $S(x^ia^j)=S(a)^iS(x)^j$ for all $0\le j,i\le p-1$.
\end{itemize}
The maps $\Delta$ and $\epsilon$ satisfy the the first relation and $\epsilon$ fulfils both relations and is therefore an algebra homomorphism. 
In the case $p=2$ this is just the Hopf algebra given in Example \ref{exa}. For $p=3$ an easy calculation shows that the map $\Delta$ satisfies the second relation $\Delta(x)^3=\Delta(x)$ and is therefore an algebra homomorphism. In this case the axioms of a Hopf algebra hold and $A_3$ is a Hopf algebra. 

\end{exa}

In order to classify third type Hopf algebras we first show that $x$ has degree $p$ and $x^p=x+\alpha_0(a^p-1)$ with $\alpha_0\in k$. To this end we will compute the coefficients of elements $x^i\otimes x^za^j$ with $0\le i,z,j\le p$ in $\Delta(x^p)$. Throughout this section $p$ is a prime number.

\begin{lemma}[Fermat's little Theorem]\label{fermat}
In a field of charactistic $p$ we have $\prod_{m=1}^{p-1}(t-m)=t^{p-1}-1$ where $t$ is a variable.
%as both are normed polynomials of degree $p-1$ with $p-1$ distinct nullstellen $\{1,\ldots,p-1\}.$
\end{lemma}
Furthermore we have to define two maps. 
\begin{lemma}\label{prime}For $j,l\in \N$ we define \[f(j,l)=\sum_{{1\le j_1\le \ldots\le j_l\le j}} \prod_{u=1}^lj_u\] and for $r,q \in \N$ with $0<q< r$\[g(r,q)=\sum_{0< j_1< \ldots <j_q <r}\prod_{u=1}^qj_u.\] Then $f(j,p-j)=g(p, p-j)=0 \mod p$ for $1 < j < p$. Furthermore $f(1,p-1)=1$ and $g(p, p-1)= (p-1)!=-1 \mod p$.
\end{lemma}
\begin{proof} 

The value of $f(j,p-j)$ is the coefficient of $t^{p-j}$ in the series expansion of
$\prod_{i=1}^{j}1/(1-it)$, where $t$ is a variable. Using Fermat's little Theorem we have
$\prod_{i=1}^{p-1}1/(1-it)=1/(1-t^{p-1})$. Thus
\[\prod_{i=1}^{j}1/(1-it)=1/(1-t^{p-1})*\prod_{i=j+1}^{p-1}(1-it)=
(\sum_{i=0}^{\infty}t^{i(p-1)})\prod_{i=j+1}^{p-1}(1-it).\] As
$\prod_{i=j+1}^{p-1}(1-it)$ is a polynomial of degree $p-j-1$ the
coefficients of $t^{p-j}$ are zero for $1<j<p$ and the coefficient
of $t^{p-1}$ is 1. 

The value of $g(p,p-j)$ is just the coefficient of $t^{j-1}$ up to sign in the
polynomial $\prod_{m=1}^{p-1}(t-m)$. Thus by Fermat's little Theorem we
have $g(p,p-j)=0 \mod p$ for $1<j<p$ and $g(p,p-1)=-1$. \end{proof}

We set $f(j,0)=g(j,0)=1$ and $f(j,l)=g(j,l)=0 $ for $l<0$ and $j>0$. 
As one can easily see the following relations are true. 
\begin{lemma}\label{fg}
We have $jf(j,l-1)+f(j-1,l)=f(j,l)$ and $(j-1)g(j-1,m-1)+g(j-1,m)=g(j,m)$ for all $m,l$ and $j>0$.
\end{lemma}
The next lemma and corollary work in an algebra $A$ over a field of characteristic $p$ generated by an element $x$ and an element $a$ such that the set $\{x^ia^j|i,j \in\N \}$ is linearly independent and $ax=xa+a^2-a$. We also fix a scalar $c(g) \in k$.  
\begin{lemma}\label{rec}
We write $(x+c(g)a)^k$ as a linear combination of elements in $\{x^ia^j| 0\le i,j \le k\}$. Let $s^{k}_{i,j}$ be the coefficient of $x^ia^j$. We define $l:=k-(i+j)$. Then we have $s^k_{i,j}= (-1)^l\binom{k}{i}f(j,l)\prod_{u=1}^j[c(g)-1+u]$ for $j>0$ and $s^k_{i,0}=\delta_{i,k}$.
%The coefficients $s^{k}_{i,j}$ of $x^ia^j$ are given by the recursion \[s^{k}_{i,j}:=s^{k-1}_{i-1,j}+(j-1+c(g))s^{k-1}_{i,j-1}-js^{k-1}_{i,j},\]where $s^1_{0,1}=c(g)$, $s^1_{1,0}=1$ and the other coefficients are zero. We define $l:=k-(i+j)$. Then we have $s^k_{i,j}= (-1)^l\binom{k}{i}f(j,l)\prod_{u=1}^j[c(g)-1+u]$ for $j>0$ and $s^k_{i,0}=\delta_{i,k}$ for $j=0$ .

%\smallskip

%Thus for a prime number $p$ in a field of characteristic $p$ we have $s^p_{p,0}=1,$ $s^p_{0,1}=c(g)$ and $s^p_{0,p}=c(g)^p-c(g)$ and all other coefficients are zero.
We write $(x\otimes a +1\otimes x)^k$ in $A\otimes_kA$ as linear combination of elements in $\{x^i\otimes x^za^j| 0\le i,j,z \le k\}$. Let $h^{k}_{i,z,j}$ be the coefficient of $x^i\otimes x^za^j$.  We define $l:=k-(z+j)$ and $m=j-i$. Then we have $h^k_{i,z,j}= (-1)^l\binom{k}{z}f(j,l)g(j,m)$ for $0<i,j$, $h^k_{0,z,j}=\delta_{j,0}\delta_{k,z}$ and $h^k_{i,z,0}=\delta_{i,0}\delta_{k,z}$. %The coefficients $h^{k}_{i,z,j}$ of $x^i\otimes x^za^j$ are given by the recursion \[h^{k}_{i,z,j}:=h^{k-1}_{i-1,z,j-1}+h^{k-1}_{i,z-1,j}+(j-1)h^{k-1}_{i,z,j-1}-jh^{k-1}_{i,z,j},\]where $h^1_{0,1,0}=h^1_{1,0,1}=1$ and zero else. We define $l:=k-(z+j)$ and $m=j-i$. Then we have $h^k_{i,z,j}= (-1)^l\binom{k}{z}f(j,l)g(j,m)$ for $0<i,j$ and $h^k_{0,z,j}=\delta_{j,0}\delta_{k,z}$ for $i=0$. In the case $j=0$ we have $h^k_{i,z,0}=\delta_{i,0}\delta_{k,z}$.Thus for a prime number $p$ in a field of characteristic $p$ we have $h^p_{0,p,0}=1$, $h^p_{p,0,p}=1$, $h^p_{1,0,p}=(p-1)!$ and $h^p_{1,0,1}=1$ and all other coefficients are zero. 

We write $(xa^{-1})^k$ as a linear combination of elements in $\{ x^ia^{-j}|0\le i,j \le k\}$. Let $c_{i,j}^k$ be the coefficient of $x^i a^{-j}$. Then we have $c_{i,j}^k= (-1)^{k-j}g(j,j-i)f(j,k-j)$ for $ 1\le i,j \le k$ and $c_{0,j}^k=c_{i,0}^k=0$.  
\end{lemma}
\begin{proof} 
By definition it is clear that $a^jx=xa^j+j(a^{j+1}-a^j)$ for all $j \ge 0$. If we set $s^k_{v,w}=0$ for $v<0 $ or $w < 0$, then the coefficients $s^{k}_{i,j}$ of $x^ia^j$ are given by the recursion 
\[s^{k}_{i,j}:=s^{k-1}_{i-1,j}+(c(g)+j-1)s^{k-1}_{i,j-1}-js^{k-1}_{i,j},\]where $s^1_{0,1}=c(g)$, $s^1_{1,0}=1$ and the other coefficients $s^1_{i,j}$ are zero. This follows directly from the fact that $x^ia^j(x+c(g)a)= x^{i+1}a^j +(c(g)+j)x^ia^{j+1}-jx^ia^j$.
First one checks that $s^k_{i,j}=0$ if $i+j>k$ by a straightforward induction using recursion. 
We will show by induction on $k$ that the stated formulae hold for the remaining case that $i+j\le k $. For $k=1$ they are true. Suppose now that the formulae are true for $k-1$. If $j=0$ we have by the recursion $s^k_{i,0}= s^{k-1}_{i-1,0}=\delta_{i-1,k-1}=\delta_{i,k}$.

We have $0\le k-(i+j)=:l$. If $j=1$ and $i< k-1$ we have \[s^k_{i,1}=(-1)^{l}\binom{k-1}{i-1}f(1,l)c(g)-(-1)^{l-1}\binom{k-1}{i}f(1,l-1)c(g)= (-1)^{l}\binom{k}{i}c(g)\] using that $f(1,r)=1$ for all $0\le r$. If $j=1$ and $i=k-1$ we have \[s^k_{k-1,1}=(k-1)f(1,0)c(g)+c(g)-(-1)f(1,-1)c(g)=kc(g)\]

If $j>1$ we have {\allowdisplaybreaks \begin{eqnarray*}
s^k_{i,j}&=&(-1)^{l}\binom{k-1}{i-1}f(j,l)\prod_{u=1}^{j}[c(g)-1+u]\\& &+(j-1+c(g))\prod_{u=1}^{j-1}[c(g)-1+u](-1)^{l}\binom{k-1}{i}f(j-1,l)\\& &-j(-1)^{l-1}\binom{k-1}{i}f(j,l-1)\prod_{u=1}^{j}[c(g)-1+u]\\&=&(-1)^{l}\binom{k-1}{i-1}f(j,l)\prod_{u=1}^{j}[c(g)-1+u]\\ & & +(-1)^{l}\binom{k-1}{i}f(j,l)\prod_{u=1}^{j}[c(g)-1+u]\\&=&(-1)^{l}\binom{k}{i}f(j,l)\prod_{u=1}^{j}[c(g)-1+u]\end{eqnarray*}} using the first identity of Lemma \ref{fg}. This proves the coefficient formula for $s$.

%Suppose now $j>0$. If we iterate by using the recursion formula we get \[s^k_{i,j}=\sum_{l=0}^{k-1}\sum_{\{r,s\in \N|r+s=k-1-l\}} s^1_{i-r,j-s}(-1)^l\binom{k-1}{r}f(j,l)\prod_{u=j-s+1}^{j}[c(g)-1+u].\] As all $s^1$ except $s^1_{0,1}$ and $s^1_{1,0}$ are zero, we only need to consider the two cases  $(r=i-1, s=j)$ and $(r=i,s=j-1)$. In both cases we have $l=k-(i+j)$. Thus \begin{eqnarray*} s^k_{i,j}&=&(-1)^l\binom{k-1}{i}f(j,l)\prod_{u=1}^{j}[c(g)-1+u]\\&\ &+(-1)^l\binom{k-1}{i-1}c(g)f(j,l)\prod_{u=2}^{j}[c(g)-1+u]\\ &=&(-1)^l\binom{k}{i}f(j,l)\prod_{u=1}^{j}[c(g)-1+u].\end{eqnarray*}If $j=0$ we deduce $s^k_{i,0}=s^1_{i-k-1,0}=\delta_{i,k}$ using iteration. This proves the first formula. 
%Thus if we take $p=k$, $j>0$ and $0<i<p$ we get $s^p_{i,j}= 0\mod p$ as $(p,i)=0 \mod p$. If $i=0$ and $j>0$ then we have $s^p_{0,j}=(-1)^{p-j}f(j,p-j)\prod_{u=1}^{j}[c(g)-1+u]$. If we use Lemma \ref{prime} we see that $s^p_{0,j}=0$ for $1<j<p$, $s^p_{0,1}= c(g)$ and $s^1_{0,p}= \prod_{u=1}^{p}[c(g)-1+u]=c(g)\prod_{u=1}^{p-1}[c(g)-u]=c(g)(c(g)^{p-1}-1)=c(g)^p-c(g)$ by Fermat's little Theorem. 

\bigskip

If we set $h^k_{r,s,t}=0$ for all $k$ and $r,0$ or $s<0$ or $t < 0$, then the coefficients $h^{k}_{i,z,j}$ of $x^i\otimes x^za^j$ are given by the recursion \[h^{k}_{i,z,j}:=h^{k-1}_{i-1,z,j-1}+h^{k-1}_{i,z-1,j}+(j-1)h^{k-1}_{i,z,j-1}-jh^{k-1}_{i,z,j},\]where $h^1_{0,1,0}=h^1_{1,0,1}=1$ and zero else. This follows from the fact that $(x^i\otimes x^za^j)(x\otimes a+1\otimes x)= x^{i+1}\otimes x^za^{j+1}+x^i\otimes x^{z+1}a^j+jx^i\otimes x^za^{j+1}-jx^i\otimes x^za^j$. By an easy induction we see that $h^k_{i,j,z}$ is zero if $j-i< 0$ or $z+j>k$. We show by induction on $k$ that the formula also holds for the other cases. We set $0\le j-i=:m$ and $0\le k-z-j=:l$.

Suppose $i=0$, then {\allowdisplaybreaks \begin{eqnarray*}h^k_{0,z,j}&=& h^{k-1}_{0,z-1,j}+(j-1)h^{k-1}_{0,z,j-1}-jh^{k-1}_{0,z,j}\\&=&\delta_{j,0}\delta_{k-1,z-1}+ (j-1)\delta_{j-1,0}\delta_{k-1,z}-j\delta_{j,0}\delta_{k-1,z}\\&=&\delta_{j,0}\delta_{k,z}.\end{eqnarray*}}Suppose that $j=0$, then $h^k_{i,z,0}= h^{k-1}_{i,z-1,0}=\delta_{i,0}\delta_{k-1,z-1}=\delta_{j,0}\delta_{k,z}$. Now we assume that $1<i\le j$, then {\allowdisplaybreaks \begin{eqnarray*} h^{k}_{i,z,j}&=&(-1)^{l}\binom{k-1}{z}f(j-1,l)g(j-1,m)+(-1)^{l}\binom{k-1}{z-1}f(j,l)g(j,m)\\ & &+(j-1)(-1)^{l}\binom{k-1}{z}f(j-1,l)g(j-1,m-1)\\& &-j(-1)^{l-1}\binom{k-1}{z}f(j,l-1)g(j,m)\\&=& (-1)^{l}\binom{k-1}{z}f(j-1,l)g(j,m)+(-1)^{l}\binom{k-1}{z-1}f(j,l)g(j,m)\\& &-j(-1)^{l-1}\binom{k-1}{z}f(j,l-1)g(j,m)\\&=&(-1)^{l}\binom{k-1}{z}f(j,l)g(j,m)+(-1)^{l}\binom{k-1}{z-1}f(j,l)g(j,m)\\&=&(-1)^{l}\binom{k}{z}f(j,l)g(j,m)\end{eqnarray*}} using the second and then the first identity of Lemma \ref{fg}.

In the case $i=1, j>1$ we have {\allowdisplaybreaks \begin{eqnarray*}h_{1,z,j}&=& \delta_{j,0}\delta_{k-1,z}
+(-1)^{l}\binom{k-1}{z-1}f(j,l)g(j,j-1)\\& &
+(j-1)(-1)^{l}\binom{k-1}{z}f(j-1,l)g(j-1,j-2)\\& &-j(-1)^{l-1}\binom{k-1}{z}f(j,l-1)g(j,j-1)\\ &=&(-1)^{l}\binom{k}{z}f(j,l)g(j,j-1),\end{eqnarray*}} if we use first that $(j-1)g(j-1,j-2)=g(j,j-1)$ and then the first identity of Lemma \ref{fg}. Finally we consider the case $i=j=1$. Then \begin{eqnarray*}h^k_{1,z,1}&=&\delta_{k,z}+(-1)^{l}\binom{k-1}{z-1}f(1,l)g(1,0)\\& &
-(-1)^{l-1}\binom{k-1}{z}f(1,l-1)g(1,0)\\&=&(-1)^{l}\binom{k}{z}\end{eqnarray*} because we have $k>z$, since otherwise $l<0$ which contradicts our assumptions. This proves the statement.

%By iterating with the recursion formula we get \[h^k_{i,z,j}=\sum_{l=0}^{k-1}\sum_{ \{s,r,m\in \N| s+r+m=k-1-l\}}h^1_{i-s, z-r,j-m-s}(-1)^l\binom{k-1}{r}f(j,l)g(j,m).\] As all $h^1$ except $h^1_{1,0,1}=h^1_{0,1,0}=1$ are zero we only need to consider the cases $(s=i, r=z-1, m=j-i)$ and $(s=i-1,r=z, m=j-i)$. In both cases we have $l= k-z-j$. Thus \begin{eqnarray*}h^k_{i,z,j}&=&(-1)^l\binom{k-1}{z-1}f(j,l)g(j,m)\\& &+ (-1)^l\binom{k-1}{z}f(j,l)g(j,m)\\&=&(-1)^l\binom{k}{z}f(j,l)g(j,m). \end{eqnarray*}Suppose $i=0$, then $h^k_{0,z,j}=h^1_{0,z-k+1,j}=\delta_{j,0}\delta_{z,k}$. Now let $j=0$. Then $h^k_{i,z,0}=h^1_{i,z-k+1,0}=\delta_{i,0}\delta_{z,k}$.
%So if we take $p=k$, $j,i>0$ and $0<z<p$ we get $h^p_{i,z,j}= 0\mod p$ as $(p,z)=0 \mod p$. If $z=0$ and $i,j>0$ we have $h^p_{i,0,j}=(-1)^{p-j}f(j,p-j)g(j,j-i)$. Lemma \ref{prime} gives us $h^p_{i,0,j}=0$ for $1<j<p$ and $h^p_{1,0,1}=f(1,p-1)g(1,0)(-1)^{p-1}=1$. By lemma \ref {prime} we have $h^p_{i,0,p}=g(p,p-i)=0 \mod p$ for $p>i>1$, $h^p_{p,0,p}=1$ and $h^p_{1,0,p}=g(p,p-1)=(p-1)!$. If $z=p$ then $h(i,p,j)=\delta_{i,0}\delta_{j,0}$ because if $j>0$ then $l$ is negative and if $i>0$ and $j=0$ then $m$ is negative. 

\bigskip

By definition $a^{-j}x=xa^{-j}+j(a^{-j}-a^{-j+1})$ for $j \in \N$. Therefore $x^ia^{-j}xa^{-1}= x^{i+1}a^{-j-1}+jx^ia^{-j-1}-jx^ia^{-j}$. So in the last case the recursion formula is given by $c_{i,j}^k=c_{i-1,j-1}^{k-1}+(j-1)c_{i,j-1}^{k-1}-jc_{i,j}^{k-1} $ with $c_{i,j}^1=\delta_{i,1}\delta_{j,1}$. Then the formula for the coefficients follows as in the previous step by induction using the identities given in Lemma \ref{fg}. 

Suppose $i=0$ or $j=0$, then by the recursion formular $c_{0,j}^k=c_{i,0}^k=0$ 
and $c_{i,j}^k=0 $ for $i>j$.  
Suppose now that $i=1$ and $j=1$. Then $c_{1,1}^k=-1c_{1,1}^{k-1}=(-1)^{k-1}$. 
Next suppose that $i=1$ and $ j>1$. Then $c_{1,j}^k= (j-1)c_{1,j-1}^{k-1}-jc_{1,j}^{k-1}=
 (j-1)(-1)^{k-j}g(j-1,j-2)f(j-1,k-j)+j(-1)^{k-j}g(j,j-1)f(j,k-j-1)=(-1)^{k-j}(j-1)![f(j-1,k-j)+j f(j,k-j-1)]=
(-1)^{k-j}g(j,j-1)f(j, k-j)$ by the second equation of \ref{fg}.
Finally we can assume $1<i \le  j \le k$. Then $c_{i,j}^k= (-1)^{k-j}[g(j-1, j-i)f(j-1,k-j)+(j-1)g(j-1,j-i-1)f(j-1,k-j)
+jg(j,j-i)f(j,k-j-1)]= (-1)^{k-j}[g(j, j-i)f(j-1,k-j)+jg(j,j-i)f(j,k-j-1)]=(-1)^{k-j}g(j, j-i)f(j,k-j)$ using first the 
identity on $g$ in \ref{fg} and then the identity on $f$ in \ref{fg}. 

\end{proof}
Using the result of the previous lemma and the same setup we can now determine the coefficients for the $p$-th power. 
\begin{cor} \label{coef}
We have $(x+c(g)a)^p=x^p+c(g)a+(c(g)^p-c(g))a^p$, $(x\otimes a +1\otimes x)^p=x^p\otimes a^p+(p-1)!x\otimes a^p+x\otimes a+1\otimes x^p$ and $(xa^{-1})^p=x^pa^{-p}-xa^{-p}+xa^{-1}$.  
\end{cor}
\begin{proof} We have to show that $s^p_{p,0}=1,$ $s^p_{0,1}=c(g)$ and
$s^p_{0,p}=c(g)^p-c(g)$ and that all other coefficients are zero by
using the formula provided in Lemma \ref{rec}. If we take $0<j,i<p$, we get $s^p_{i,j}= 0\mod p$ as $\binom{p}{i}=0 \mod p$.
If $i=0$ and $0<j\le p$, then we have
$s^p_{0,j}=(-1)^{p-j}f(j,p-j)\prod_{u=1}^{j}[c(g)-1+u]$. Lemma
\ref{prime} proves that $s^p_{0,j}=0$ for $1<j<p$ and $s^p_{0,1}= c(g)$.
By Fermat's little Theorem we get $s_{0,p}^p=
\prod_{u=1}^{p}[c(g)-1+u]=c(g)\prod_{u=1}^{p-1}[c(g)-u]=c(g)(c(g)^{p-1}-1)=c(g)^p-c(g)$.
Furthermore we have $s_{p,j}^p=\delta_{j,0}$ since if $j>0$ then $l<0$
and $s_{i,0}^p=\delta_{i,p}$ by Lemma \ref{rec}. Finally
$s_{i,p}^p=0$ for $i>0$ since then $l<0$. This proves the first identity.

\bigskip 

In order to prove the second identity, we have to show that $h^p_{0,p,0}=1$, $h^p_{p,0,p}=1$, $h^p_{1,0,p}=(p-1)!$ and $h^p_{1,0,1}=1$ and all other coefficients are zero. We first consider the case $j,i>0$. If $0<z<p$ we get $h^p_{i,z,j}= 0\mod p$ as $\binom{p}{z}=0 \mod p$. 
If $z=0$ and $i,j>0$ we have
$h^p_{i,0,j}=(-1)^{p-j}f(j,p-j)g(j,j-i)$. Lemma \ref{prime} gives us
$h^p_{i,0,j}=0$ for $1<j<p$ and $h^p_{1,0,1}=(-1)^{p-1}
f(1,p-1)g(1,0)=1$. By Lemma \ref {prime} we have
$h^p_{i,0,p}=g(p,p-i)=0 \mod p$ for $1<i<p$, $h^p_{p,0,p}=1$ and $h^p_{1,0,p}=g(p,p-1)=(p-1)!$. If $z=p$ then $h(i,p,j)=0$ because $l<0$. 

The values of the coefficients for $i=0$ or $j=0$ follow immediately from \ref{rec}.

\bigskip

By \ref{rec} we have $c_{0,j}^p=c_{i,0}^p=0$ and $c_{i,j}^p=0$ for all $1< i < p$ using \ref{prime}. Furthermore $c_{i,1}^p=\delta_{1,i}$ as $g(1,1-i)=0$ for $i > 1$ and $c_{i,p}^p=0 $ for $1<i <p$, $c_{1,p}^p=-1 $ by \ref{prime} and $c_{p,p}^p=1$. Therefore the equation holds.

\end{proof}
One can easily see that \ref{rec} and \ref{coef} can be applied to our Hopf algebra $H$ in order to deduce some necessary conditions for the degree $n$ of $x$. 

\smallskip

Similarly to \cite[Lemma 1]{KR} we prove the following property
of the degree $n$ of $x$. 
If we set $V_0=H_0$ and $V_k:=\sum_{t=0}^{k}H_0x^t$, then all conditions of
remark \ref{free} are satisfied and $\rho_k(x^k+V_{k-1})=1\otimes
(x^k+V_{k-1})$  as $h_{j,z,i}^k=\delta_{i,0}\delta_{j,0}\delta_{k,z}$ for $z \ge k$ and $x^t H_0 \subset V_t$ by \ref{rec}. We have $Q_k\not=\{0\}$ for all $k< n$ and $Q_n=\{0\}$. As $x^{k}+V_{k-1}$ forms
a basis for $Q_k$ as $H_0$-module, $n$ is the smallest integer such that $1,x, \ldots ,x^{n-1}$ form a basis for $H$ over $H_0$.   

\begin{lemma}\label{cond} 
Let $H$ be a pointed rank one Hopf algebra of third type.
Let $n$ be the degree of $x$ and $x^n=\sum_{i=0}^{n-1} \alpha_i x^i$
with $\alpha_i \in kG$. Then the coefficients
$\alpha_i$ with $1\le i$ are in $k$ and $\alpha_0=\alpha(a^n-1)$ for some $\alpha\in k$. Furthermore $n=p^q$ for a positive integer $q$ .
\end{lemma}
\begin{proof}We use that $x^k\otimes x^j$ with $0\le k,j\le n-1$ is 
free over $kG\otimes_k kG$ and that the set $\{x^kg|g \in G,\ 0\le k\le n-1\}$ is linearly independent. We have $(x\otimes a+1\otimes x)^n=\Delta(x)^n= \sum_{i=0}^{n-1}\Delta(\alpha_i) \Delta (x)^i$. If we compare the coefficients in $H_0\otimes H_0$ of $1\otimes x^i $ for all $0\le i \le n-1$ on both sides this gives us $ 1 \otimes \alpha_i = \Delta(\alpha_i)$ for $1\le i \le n-1$ and $\alpha_0\otimes a^n+1 \otimes \alpha_0 =\Delta(\alpha_0)$. Thus $\alpha_i \in k$ for $1\le i \le n-1$ and $\alpha_0=\alpha(a^n-1)$ for an $\alpha \in k$.
 For the second part we compare the coefficients of $x^r\otimes x^{n-r}a^r$ for $1\le r
\le n-1$ using the
coefficients computed in Lemma \ref{rec}. We have $h^n_{r,n-r,r}=\binom{n}{r}$
on the left hand side and $h^i_{r,n-r,r}=0$ for all $1\le i\le n-1$ on
the right hand side. Thus $\binom{n}{r}=0\mod p $ for all $1\le r \leq n-1$ and $n$ has to be a $p$-power. 

\end{proof}
 Now we show that $n=p$ using the results for the coefficients of the term $\Delta(x)^p$ and the previous lemma. 

\begin{lemma}\label {third}Let $H$ be a pointed rank one algebra of third type. Then $x^p\in H_1$ and $p$ is the degree of $x$. We have $x^p=x+\alpha(a^p-1)$ for some $\alpha\in k$.
\end{lemma}
\begin{proof}
We have $\Delta(x^p)=x^p\otimes a^p+(p-1)!x\otimes a^p+x\otimes a+1\otimes x^p$ from Corollary \ref{coef}. Thus $x^p \in H_1$ by the definition of $H_1$. Using Lemma \ref{cond} we know that the degree of $x$ is $p$.
By part two of Lemma \ref{cond} we have $x^p=\alpha_1x+\alpha(a^p-1)$ with
$\alpha_1,\alpha \in k$. Comparing the coefficients
of $x\otimes a$ in $\Delta(x)^p $ with the coefficient in
$\Delta (\alpha_1x+\alpha(a^p-1))$ gives $\alpha_1=1$ if $a^p\not = a$. By \ref{types} we know that the order of $a$ is divisible by $p$ and therefore $a^p\not= a$. 
\end{proof}If $a$ has order $p$, then $\alpha=0$. If $k$ is
algebraically closed, then we can choose $\beta \in k$ such that
$\beta-\beta^p=\alpha$ and set $\bar x=x+\beta(a-1)$. In this case $\bar x\in P_{a,1}\setminus H_0$ and $g\bar x g^{-1}= \bar x +c(g)(a-1)$ as $a$ is central. Furthermore $\bar x^p= (x+\beta a)^p-\beta^p= x^p+\beta a+(\beta^p-\beta)a^p-\beta^p= \bar x$. We can therefore assume without loss of generality that $\alpha=0$.

In the next definition we introduce Hopf algebras parametrized by a tuple $E$.
\begin{defi}[third type]\label{H_E} Let $E:=(G,a,c, \alpha)$ be a tuple, where $G$ is a finite group,
$a\in Z(G)$ and $a$ has an order divisible by $p$ and $c:G\to (k,+)$ is
a group homomorphism. The Hopf algebra $H_E$ is the algebra with basis
$(gx^i)_{\{g\in G, 0\le i\le p-1\}}$ and relations
$gxg^{-1}=x+c(g)(a-1)$ and $x^p=x+\alpha(a^p-1)$.
The Hopf algebra structure is determined by the linear maps $\Delta$, $\epsilon $ and $S$ with $\Delta(gx^i)=(g\otimes g)(x\otimes a+1\otimes x)^i$, $\epsilon(gx^i)=\delta_{i,0}$ and $S(gx^i)=(-xa^{-1})^ig^{-1}$ for all $g\in G$ and $0\le i \le p-1$. If $a^p \not =1$ we assume that $c(g)\in F_p$ for all $g\in G$.  
\end{defi}
The proof of the next lemma shows that $H_E$ with $a^p\not =1$ is well-defined if and only if $c(g)\in F_p$ for all $g\in G$.
\begin{lemma} The Hopf algebra $H_E$ is well-defined. 
\end{lemma}
\begin{proof}
Using the results of Corollary \ref{coef}, we get the following equation
\begin{eqnarray*} x^p+c(g)(a-1)&=&gx^pg^{-1}=(gxg^{-1})^p\\
&=&(x+c(g)a-c(g))^p=(x+c(g)a)^p-c(g)^p\\
&=&x^p+c(g)a+(c(g)^p-c(g))a^p-c(g)^p.
\end{eqnarray*}
If $a^p=1$, then the equation holds for any value of $c(g)$. If $a^p\not =
1$, then $gx^pg^{-1}=(gxg^{-1})^p$ is true if and only if $c(g)^p=c(g)$,
which is the case if and only if $c(g)\in F_p$ for all $g\in G$. Then $H_E$ is well-defined as an algebra. 
Using Theorem \ref{third} and the expression for $\Delta(x)^p$ in Corollary \ref{coef}, we see that $\Delta$ and $\epsilon$ are algebra homomorphism. 

Let $S$ be a linear map with $S(gx^i)=(-xa^{-1})^ig^{-1}$  on the basis $\{gx^i|g\in G, 0\le i \le p-1\} $. By \ref{coef} and \ref{third} we have $S(x)^p=S(x)+\alpha(a^{-p}-1)$ and $S(gxg^{-1})=S(g^{-1})S(x) S(g)$ for all $g\in G$. Therefore $S$ is an anti-algebra homomorphism. It is easy to prove that the axiom is satisfied for all elements of this basis. 
%Therefore note that the following general property holds: Let $z=\sum_{k=1}^t a_k\otimes b_k \in H \otimes H$ and let $f:H \to H$ be an anti-algebra homomorphism. Suppose $m [ (f\otimes \id )(z)]=0$. Then $m [(f\otimes \id)(z^2)]=\sum_{k=1}^t\sum _{l=1}^tf(a_la_k)b_lb_k=\sum_{k=1}^t\sum _{l=1}^tf(a_k)f(a_l)b_lb_k=\sum_{k=1}^t f(a_k)[m(f \otimes \id )(z)]b_k=0$. We have $\Delta(x^i)=\Delta(x)^i$. Therefore the previous remark and induction prove that $m[(S\otimes \id) (\Delta(x^i)]=0$ for all $1\le i \le p-1$. Let $\Delta(x^i)=x_{(1)}\otimes x_{(2)}$ for $1\le i \le p-1$. Then $m[(S\otimes \id) (\Delta(gx^i))]=S(gx_{(1)})gx_{(2)}=S(x_{(1)})g^{-1}gx_{(2)}=S(x_{(1)})x_{(2)}=0$ for all $ g\in G$. Therefore $S$ is the antipode of $H_E$. 
\end{proof}We can now classify Hopf algebras of third type. 
\begin{theo}Let $H$ be a third type, then $H$ is isomorphic to a Hopf
algebra $H_E$ for a tuple $E= (G(H),a,c, \alpha)$.
\end{theo}
\begin{proof} 
This follows immediately from Lemma \ref{third}, Definition \ref{H_E} and Corollary \ref{c}. \end{proof}
The converse is also true. 
\begin{lemma}The Hopf algebras $H_E$ are pointed, rank one and of third type.
\end{lemma}
\begin{proof}Exactly as in \ref{rankone} we can see that $H_E$ is pointed,
as it is generated by skew primitive elements.
We want to show that it is rank one.
Let $z\in H_E\setminus kG$ be an element with $\Delta(z)=z\otimes g +1\otimes z $ for some $g\in G$.
Then $z$ has a unique presentation as $z=\sum_{t=0}^{p-1}x^tb_t$ with $b_t\in kG$ for all $0\le t \le p-1$.
These conditions give the following equation 
\begin{eqnarray*} \sum_{t=0}^{p-1}\sum_{i,z,j=0}^th^t_{i,z,j}(x^i\otimes
x^z a^j)\Delta(b_t)&=&\sum_{t=0}^{p-1}\Delta(x)^t\Delta(b_t)\\&=& \Delta(z)=z\otimes g +1\otimes z\\&=&\sum_{t=0}^{p-1} (x^t\otimes 1)(b_t\otimes g)+(1\otimes x^t)(1\otimes b_t).
\end{eqnarray*}
We compare the coefficients of $1\otimes x^z$ for $0\le z \le p-1$ on both sides. We have $h^t_{0,z,j}=\delta_{j,0}\delta_{t,z}$ for all $0\le t\le p-1$. Therefore $\Delta(b_z)=1\otimes b_z$ for $1\le z\le p-1$. This forces $b_z\in k$ for all $1\le z \le p-1$. In the case $z=0$ we have $\Delta(b_0)=b_0\otimes g+1\otimes b_0$ and therefore $b_0=b(g-1)$ for some $b\in k$. Let $n$ be maximal with $b_n \not =0 $.
We compare now the coefficients of $x^n\otimes 1$. We have $h^t_{n,0,j}=\delta_{j,n}\delta_{t,n}$ for $t\le n$ by \ref{rec}. This
gives us $(1\otimes a^n)\Delta(b_n)=b_n\otimes g$. Therefore $g=a^n$. Suppose $n>1$.
The coefficient of $x\otimes x^{n-1}a$ on the left hand side is $\Delta(b_n)n \not =0$ as $h^t_{1,n-1,1}=\delta_{n,t} \binom{n}{n-1}$ for $t\le n$.
As $x\otimes x^{n-1}a$ does not appear on the right hand side, this is a
contradiction. Thus $z=b_1 x+b(a-1)$ for $b_1 \in k^*$ and $b \in k$.

For all $h\in P_{v,u}\setminus H_0$ with $v,u \in G$ we have $h=0$ for $u^{-1}v \not= a$ and else $h=u(b_1 x+b(a-1))$. As $H_1= H_0+\sum_{v,u \in G(H)}P_{v,u}\subset H_0+H_0x$ and $x\not \in H_0$,
the dimension of $k\otimes_{H_0}H_1$ is two. Thus $H_E$ is a rank one Hopf algebra.
It is clearly of third type as the skew point $a$ does not have an eigenvector in $P_{a,1} \setminus H_0$ of the conjugation action. 
\end{proof}
In the following Lemma we determine precisely under which conditions two third type Hopf algebras presented by different tuples are isomorphic. 
\begin{lemma}Let $H_E$ and $H_{E'}$ be two Hopf algebras of third type
over an algebraically closed field given by tuples $E$ and $E'$. Then they are isomorphic if and only if
there is a group isomorphism from $f:G \to G'$ such that $f(a)=a'$ and
$c'\circ f=c$.
\end{lemma}
\begin{proof}Let $f:G\to G'$ be a group isomorphism with $f(a)=a'$ and $c'\circ f=c$. We define the linear map $F:H_E\to H_{E'}, gx^i \to f(g)x'^i$ which is bijective. 
 This gives an Hopf algebra isomorphism.  

Conversely, we suppose now that $H_E$ and $H_{E'}$ are isomorphic as
Hopf algebras via a map $F$. As group-like elements are mapped to
group-like elements, the map $f:=F|_{G}:G\to G'$ is a group isomorphism. Since elements in $H_1\setminus H_0$ are mapped to elements in $H_1'\setminus H_0'$ and skew
primitives to skew primitives, we have $f(a)=a'$ and $F(x)=\gamma
x'+\beta(a'-1)$ for some $\beta, \gamma \in k$. We have $ \gamma
x'+(\beta+\gamma)(a'-1)=a'F(x)a'^{-1}=F(axa^{-1})=\gamma
x'+(\beta+1)(a'-1)$. Thus $\gamma=1$.
We have $x'+c'(f(g))(a'-1)+\beta(a'-1)=f(g)F(x)f(g)^{-1}=F(gxg^{-1})=x'+c(g)(a'-1)+\beta(a'-1)$.
As $a'\not =1$ we have $c(g)=c'(f(g))$ for all $g\in G$. 
\end{proof}  
In order to analyse the representation theory of third type
Hopf algebras we first analyse its subalgebra $k[x]/\langle x^p-x
\rangle$. 
\begin{lemma}\label{kx}Let $k[x]_p:=k[x]/\langle x^p-x \rangle$ be a polynomial ring with indeterminate $x$ and relation $x^p=x$. Let $k$ be a field of characteristic $p$. Then $k[x]_p$ viewed as a left $k[x]_p$-module is the direct sum of $p$ simple one dimensional modules $S_c$ on which $x$ acts as scalar $c$ with $c \in F_p$.  

%Then $k[x]_p$ has central orthogonal primitive idempotents $e_c=-\sum_{i=1}^{p-1}c^ix^i$ for all $c\in F^*_p$ and $e_0=1-x^{p-1}$ such that $\sum_{c\in F_p}e_c=1$. Each idempotent induces a simple one dimensional module $S_c$ on which $x$ acts as scalar $c$.  
\end{lemma}
\begin{proof}We have $x^p-x=\prod_{c\in F_p}(x-c)$. Therefore $k[x]/\langle x^p-x \rangle $ is isomorphic to $\prod_{x\in F_p}k[x]/\langle x-c \rangle$ by the Chinese Reminder Theorem. Notice that $k[x]/\langle x-c \rangle$ is isomorphic to $k$ for any $c\in F_p$. Therefore $k[x]_p$ is semi-simple and the simple modules are as stated.  

%We have $c^jx^j e_c= e_c$ for all $c \in F^*_p$. Thus $e_c^2=e_c$. We have $(x^{p-1})^2=x^{p}\cdot x^{p-2}=x^{p-1}$. Thus $(1-x^{p-1})^2=1-x^{p-1}$. As $k[x]_p$ is commutative all idempotents are central. We have $\sum_{c\in F^*_p}e_c=-\sum_{i=1}^{p-1}\sum_{c\in F^*_p}c^i x^i= -(p-1)x^{p-1}=x^{p-1}$ as follows:The group $F^*_p$ is cyclic of order $p-1$ generated by some $y$. Thus we have $\sum_{c\in F^*_p}c^{p-1}=\sum_{c\in F^*_p}1=p-1$ and $\sum_{c\in F^*_p}c=p(p-1)/2=0 \mod p$. Let $1<i<p-1$ and $d=\mbox{gcd}(p-1,i)$ then $ \sum_{c\in F^*_p}c^i=\sum_{j=1}^{p-1}y^{ij}=d\sum_{j=1}^{p-1/d}y^{dj}$. We set $q:=\sum_{j=1}^{p-1/d}y^{dj}$, then $q(1-y^d)/(1-y)=\sum_{i=0}^{d-1}y^{i}q=\sum_{i=1}^{p-1}y^i=p(p-1)/2=0 \mod p$. Therefore $q=0 \mod p$. 

%\smallskip

%This gives us $\sum_{c\in F_p}e_c=1$ and the idempotents are therefore orthogonal. As there are $p$ different orthogonal idempotents the modules $S_c:=k[x]_pe_c$ are all one dimensional and therefore simple. The action of $x$ on $S_c$ is given by multiplication with the scalar $c$. Finally $k[x]_p$ is semi-simple as it is the direct sum of $p$ simple modules.  
\end{proof}
Let $\{e_c|c\in F_p \}$ be the set of orthogonal central primitive idempotents in $k[x]_p$ such that $k[x]_p\cong \bigoplus_{c\in F_p} k[x]_p e_c$ and $ k[x]_p e_c \cong S_c$.  
We can now describe the projective and simple modules of certain third type Hopf algebras.
\begin{lemma}\label{subhopf}Let $A$ be of third type, with $G=\langle
a\rangle $ cyclic of $p$-power order and $\alpha=0$. Then $A$ is indecomposable and its simple modules are one-dimensional.
There are $p$ simple modules $T_c$ indexed by $c\in F_p$ on which $x$ acts by scalar multiplication with $c$ and $a $ as the identity. The projective indecomposable modules of $A$ are $Ae_c$. 
\end{lemma}
\begin{proof}We have $A=\oplus_{c\in F_p} Ae_c$. The submodules $Ae_c$ are clearly indecomposable as left $k\langle a \rangle$-modules.
We set $y:=1-a$. Then a composition series of $Ae_c$ is
${P_c}^{p^r-1}\subset {P_c}^{p^r-2} \subset \cdots \subset {P_c}^0$
where ${P_c}^i:=\langle y^je_c| i\le j \le p^r-1 \rangle$.
 The quotients are one dimensional and $a$ acts on them as the identity.
We have $xy^ie_c=y^ixe_c+i(y^{i+1}-y^i)e_c$ and therefore $x$ acts on ${P_c}^i/ {P_c}^{i+1}$ by skalar multiplication with $c-i$.
As any two projective indecomposable modules have a common composition factor, the algebra $A$ is indecomposable. 
\end{proof}

\bigskip

 If $a$ has order $p$, the algebra in \ref{subhopf} is
isomorphic to a restricted $p$-Lie algebra constructed in the following
way: We choose $L:=\langle u,v \rangle$ and $[u,v]=-v$ and define the
$p$-map $u^{[p]}:=u$ and $v^{[p]}:=0$. It is easy to see that
$(\mbox{ad}u)^p=\mbox{ad}u^{[p]}=\mbox{ad}u$ and
$(\mbox{ad}v)^p=\mbox{ad}v^{[p]}=0$. In $U^{[p]}(L)$ we
have $u^p=u$ and $v^p=0$. If we set $x:=u$ and $a^{-1}:=1-v$ we have
$a^p=1$ and $axa^{-1}=x+(a-1)$. Therefore $U^{[p]}(L)$ is isomorphic to $A$. 

\bigskip

In general we can describe third type Hopf algebras as follows.

\bigskip 

Let $H$ be a third type Hopf algebra. The group $G$ acts on
the two dimensional subspace $A$ spanned by $a$ and $x$ via conjugation.
This gives by \ref{c} the following representation $\rho:G\to Gl_n(A), g \mapsto \bigl ( \begin{smallmatrix} 1 & c(g) \\ 0 & 1 \end{smallmatrix} \bigr ) .$

Let $Z:=\{g\in G|c(g)=0\}$ be the kernel of $\rho$. The image of $c$ is isomorphic to a finite additive subgroup of $k$. Furthermore all elements of $G/Z$ have order $p$. Therefore $G/Z$ is a elementary abelian $p$-group.
If $a$ does not have order $p$, then $|\Im c|=p$ and therefore $G/Z$ is
cyclic of order $p$ and generated by $aZ$. Thus $G\cong Z\times_{\langle
a^{p}\rangle} \langle a \rangle\cong Z\times \langle a \rangle/\langle b\times a^{p} \rangle$ with $b=a^{-p}\in Z$.
Let $\{e_i|1\le i\le s\}$ be the set of orthogonal central primitive idempotents in $kZ$ and $B_i:=kZe_i$ its blocks. Then $e_i$ commutes with all elements of $H$. Therefore $H=\oplus_{i=1}^s He_i$ is a decomposition into subalgebras $H_i:=He_i$.  Set $\bar a:= ae_i$, $\bar x:=xe_i$ and let $A_i$ denote the subalgebra generated by $\bar a$ and $\bar x$. Then $H_i\cong A_i\otimes_kB_i/\langle e_i\otimes_k {\bar a}^p- {\bar a}^p\otimes_k e_i \rangle $. 
If ${\bar a}^p=e_i$, then the algebra $A_i$ is
isomorphic to the indecomposable algebra in \ref{subhopf}. In this case $H_i\cong
A_i\otimes_k B_i$. The algebra structure of $B_i$ can be arbitrarily complicated.

By an argument similar to the argument before \ref{nil rep} the algebra $H_i$ is indecomposable. 

% As $Z(A_i)=k\langle \bar a^p\rangle \subset Z(B_i)$, $Z(A_i)$ is indecomposable. Therefore $A_i$ is indecomposable as well and $\bar a^p$ is a scalar $we_i$ with $w\in k$  or $Z(A_i)\cong k C_{p^t}$. As the tensor product of two indecomposable algebras is indecomposable, $A_i\otimes B_i$ is indecomposable and by a lifting argument $H_i$ is indecomposable as well.  

\end{document}